# Homing missile guidance law with imperfect measurements and imperfect information about the system.

Jaykov Foukzon

Israel Institute of Technology

A. A. Potapov

IRE Russian Academy of Sciences

**Abstract.** In this study, the generic imperfect dynamic models of air-to-surface missiles are given in addition to the related simple guidance law.

## 1. Mathematical Challenge : Creating a Game Theory that Scales.

*What new scalable mathematics is needed to replace the traditional Partial Differential Equations (PDE) approach to differential games?*

Many stochastic optimal control problems essentially come down to constructing a function $u(t,x)$ that has the properties:

(**1**)

$$u(t,x) = \inf_\alpha \mathbf{E}[\bar{\mathbf{J}}(\mathbf{X}^x_{t,D}(\omega), \alpha(t))], \qquad (1.1)$$

(**2**)

$$u(t,x) = \inf_\alpha \mathbf{E}\left[ \bar{\mathbf{J}}\left( \{\mathbf{X}^x_{s,D}(\omega)\}_{s\in[0,t]}, \{\alpha(s)\}_{s\in[0,t]} \right) + u(t, \mathbf{X}^x_{t',D}(\omega)) \right], \qquad (1.2)$$

$$t < t', \alpha(t) \in U, U \subsetneq \mathbb{R}^n,$$

where $\bar{\mathbf{J}}$ is the termination payoff functional, $\alpha(t)$ is a control and $(\mathbf{X}^x_{t,D})_{t\geq 0}$ is an Markov process governed by some stochastic Ito's equation driven by a white noise of the form

(**3**)

$$\mathbf{X}_{t,D}^{x}(\omega) = x + \int_{0}^{t} g(\mathbf{X}_{s,D}^{x}(\omega), \boldsymbol{\alpha}(t))ds + \sqrt{D}\mathbf{W}(t,\omega). \tag{1.3}$$

Traditionally, the function $u(t,x)$ has been computed by way of solving the associated *Bellman equation*, for which various numerical techniques mostly variations of the finite difference scheme have been developed. Another approach, which takes advantage of the recent developments in computing technology and allows one to construct the function $u(x,t)$ by way of backward induction governed by Bellman's principle such that described in [1]. In paper [1] equation (1.3) is *approximated* by an equation with affine coefficients which admits an explicit solution in terms of integrals of the exponential Brownian motion. In approach proposed in paper [2] we have replaced Eq.(1.3) by Colombeau-Ito's equation (1.4)

$$(\mathbf{X}_{t,D,\epsilon}^{x,\varepsilon}(\omega,\varpi))_{\epsilon} = x +$$
$$\left(\int_{0}^{t} g_{\epsilon}(\mathbf{X}_{s,D,\epsilon}^{x,\varepsilon}(\omega,\varpi), \boldsymbol{\alpha}(t))ds\right)_{\epsilon} + \sqrt{D}\int_{0}^{t}(\mathbf{w}_{\epsilon}(s,\varpi)ds)_{\epsilon} + \sqrt{\varepsilon}\mathbf{W}(t,\omega), \epsilon \in (0,1], \tag{1.4}$$
$$\varepsilon \ll 1.$$

where $\mathbf{w}_{\epsilon}(t,\varpi)$ is a regularized white noise.
Fortunately in contrast with Eq.(1.3) one can solve Eq.(1.4) **without any approximation** using strong large deviations principle [2]. In this paper (*for shortness only*) we considered quasi stochastic case, i.e. case $D = 0$.
**Statement of the novelty and uniqueness of the proposed idea:** A new approach, which is proposed in this paper allows one to construct the Bellman function $\mathbf{V(t,x)}$ and optimal control $\boldsymbol{\alpha}(t,\mathbf{x})$ directly, i.e., without any reference to the Bellman equation, by way of using strong large deviations principle for the solutions Colombeau-Ito's SDE [2].
**2. Proposed Approach.** Let us consider an $m$-persons Colombeau differential game $\mathbf{CDG}_{m;T}(\mathbf{f,g,y})$ with nonlinear dynamics:

$$(\dot{x}_{\epsilon}(t))_{\epsilon} = \mathbf{f}_{\epsilon}(t, (x_{\epsilon}(t))_{\epsilon}, \boldsymbol{\alpha}(t)); \forall t : (x_{\epsilon}(t))_{\epsilon} \in \widetilde{\mathbb{R}}^{n},$$

$$x(0) = x_0, t \in [0,T], \epsilon \in (0,1], \boldsymbol{\alpha}(t) = (\alpha_1(t),\ldots,\alpha_m(t)), \alpha_i(t) \in U_i \subsetneq \mathbb{R}^{k_i}, \tag{2.1}$$
$$i = 1,\ldots,m,$$

and $m$-persons Colombeau differential game $\mathbf{CIDG}_{m;T}(\mathbf{f,g,y,w_1,w_2})$ with imperfect information [5],[6]:

$$(\dot{x}_\epsilon(t))_\epsilon = \mathbf{f}_\epsilon(t,(x_\epsilon(t)+\delta(t))_\epsilon,\alpha(t,\beta(t))); \forall t : (x_\epsilon(t))_\epsilon \in \widetilde{\mathbb{R}}^n, \qquad (2.1')$$

Here $\widetilde{\mathbb{R}}$ is a ring of Colombeau's generalized numbers [10-11], $\widetilde{\mathbb{R}}^n = \widetilde{\mathbb{R}} \times \ldots \times \widetilde{\mathbb{R}}; t \to \alpha_i(t)$ is the control chosen by the $i$-th player, within a set of admissible control values $U_i$, and the payoff for the $i$-th player is:

$$(\bar{\mathbf{J}}_{\epsilon,i})_\epsilon = \left( \int_0^T g_{\epsilon,i}(t, x_{\epsilon,1}(t), \ldots, x_{\epsilon,n}(t); \alpha_1(t), \ldots, \alpha_m(t)) dt \right)_\epsilon + \sum_{i=1}^n [(x_{\epsilon,i}(T))_\epsilon - y_i]^2. \quad (2.2)$$

where $t \mapsto ((x_{\epsilon,1}(t))_\epsilon, \ldots, (x_{\epsilon,n}(t))_\epsilon)$ is the trajectory of the Eq.(1.1). Optimal control problem for the $i$-th player is :

$$(\bar{\mathbf{J}}_{\epsilon,i})_\epsilon = \left( \min_{\alpha_i(t)} \left( \max_{\alpha_j(t), j \neq i} \mathbf{J}_i^\epsilon \right) \right)_\epsilon. \qquad (2.3)$$

Let us consider now a family $(\mathbf{X}_t^{\varepsilon,\epsilon})_\epsilon$ of the solutions Colombeau-Ito SDE (CSDE):

$$d(\mathbf{X}_t^{\varepsilon,\epsilon})_\epsilon = (\mathbf{b}_\epsilon((\mathbf{X}_t^{\varepsilon,\epsilon})_\epsilon, t))_\epsilon + \sqrt{\varepsilon} d\mathbf{W}(t); (\mathbf{X}_0^{\varepsilon,\epsilon})_\epsilon = x_0 \in \widetilde{\mathbb{R}}^n,$$

$$(2.4)$$

$$t \in [0,T], \epsilon \in (0,1], \varepsilon \ll 1,$$

where $\mathbf{W}(t)$ is $n$-dimensional Brownian motion, $\mathbf{b}_\epsilon(\circ, t) : \mathbb{R}^n \to \mathbb{R}^n$ is a polinomial transform,i.e.
$b_{\epsilon,i}(x,t) = \sum_{|\alpha|} b_{\epsilon,\alpha}^i(x,t) x^\alpha, \alpha = (i_1, \ldots, i_k), |\alpha| = \sum_{j=0}^k i_j, i = 1, \ldots, n.$  **Definition1.**
CSDE (1.4) is $\widetilde{\mathbb{R}}$-dissipative if there is exist Lyapunov function $V(x,t)$ and finite Colombeau constant $\widetilde{\mathbb{R}} \ni [(C_\epsilon)_\epsilon] = \widetilde{\mathbf{C}} > 0$ such that:
**(1)** $\forall \epsilon [\epsilon \in [0,1]], \quad \forall x, \|x\| \geq \widetilde{\mathbf{R}} : [\dot{V}(x,t; \mathbf{b}_\epsilon)] \leq -\widetilde{\mathbf{C}} \cdot V(x,t),$ where

$$(\dot{V}(x,t;\mathbf{b}_\epsilon))_\epsilon \triangleq \frac{\partial V(x,t)}{\partial t} + \sum_{i=1}^{n} \frac{\partial V(x,t)}{\partial x_i}(b_{i,\epsilon}(x))_\epsilon.$$

(2)

$$V_r(x,t) = \lim_{r \to \infty} \left( \inf_{\|x\|>r} V(x,t) \right) = \infty.$$

**Theorem.1.Main result.**( "**Strong large deviations principle**").[2] For any solution $\mathbf{X}_t^\varepsilon = (X_{1,t}^\varepsilon, \ldots, X_{n,t}^\varepsilon)$ dissipative CSDE (2.4) and $\mathbb{R}$ valued parameters $\lambda_1, \ldots, \lambda_n$ there exists constant $C' \in \mathbb{R}_+$, such that:

$$\liminf_{\varepsilon \to 0} \mathbf{E}\left[ \|(\mathbf{X}_t^{\varepsilon,\epsilon})_\epsilon - \lambda\|^2 \right] \leq C' \|\mathbf{U}(t,\lambda)\|^2. \tag{2.5}$$

Where $\mathbf{U}(t,\lambda) = (U_1(t,\lambda), \ldots, U_n(t,\lambda)), \lambda = (\lambda_1, \ldots, \lambda_n)$ the solution of the *linear master equation*:

$$\dot{\mathbf{U}}(t,\lambda) = \mathbf{J}[\mathbf{b}(\lambda,t)]\mathbf{U} + \mathbf{b}(\lambda,t), \mathbf{U}(0,\lambda) = x_0 - \lambda, \tag{2.6}$$

where $\mathbf{J} = \mathbf{J}[\mathbf{b}(\lambda,t)]$ is Jacobian, i.e. $\mathbf{J}$ is a $n \times n$ -matrix: $\mathbf{J}[\mathbf{b}(\lambda,t)] = [\partial b_i(x,t)/\partial x_j]|_{x=\lambda}$.

**Remark.1.** We note that $\forall \epsilon \in (0,1]:$ $(\delta_\epsilon(t))_\epsilon \triangleq \liminf_{\varepsilon \to 0} \mathbf{E}\left\| (\mathbf{X}_t^{\varepsilon,\epsilon})_\epsilon - (\mathbf{X}_t^{0,\epsilon})_\epsilon \right\| \neq 0.$

**Example.1.Nonlinear Ito SDE with a 'small' white noise.**

$$\dot{\mathbf{X}}_t^\varepsilon = -a(\mathbf{X}_t^\varepsilon)^3 - b(\mathbf{X}_t^\varepsilon)^2 - c\mathbf{X}_t^\varepsilon - \sigma \cdot t^n - \varkappa \cdot t^m \cdot \sin(\Omega t^k) + \sqrt{\varepsilon}\dot{\mathbf{W}}(t),$$

$$x(0) = x_0, t \in [0, T], \varepsilon \ll 1, 0 < a; \tag{2.6}$$

$$\dot{\mathbf{X}}_t^0 = -a(\mathbf{X}_t^0)^3 - b(\mathbf{X}_t^0)^2 - c\mathbf{X}_t^0 - \sigma \cdot t^n - \varkappa \cdot t^m \cdot \sin(\Omega t^k).$$

From a general master equation (2.6) one obtain the next linear master equation:

$$\dot{u}(t) = -(3a\lambda^2 - 2b\lambda - c)u(t) - (a\lambda^3 - b\lambda^2 - c\lambda) - \sigma \cdot t^n - \varkappa \cdot t^m \cdot \sin(\Omega t^k),$$

$$u(0) = x_0 - \lambda. \tag{2.7}$$

From the differential master equation (2.7) one obtain *transcendental master equation*

$$(x_0 - \lambda(t))\exp[-(3a\lambda^2(t) + 2b\lambda(t)) \cdot t] -$$

$$-\int_0^t (\sigma\tau^n + \varkappa\tau^m \sin(\Omega \cdot \tau) + a\lambda^3(t) + b\lambda^2(t))\exp[-(3a\lambda^2(t) + 2b\lambda(t))(t - \tau)]d\tau = 0. \tag{2.8}$$

**Numerical simulation:**

$$\dot{\mathbf{X}}_t^0 = -a(\mathbf{X}_t^0)^3 - b(\mathbf{X}_t^0)^2 - c\mathbf{X}_t^0 - \sigma \cdot t^n - \varkappa \cdot t^m \cdot \sin(\Omega t^k), X(t) = \mathbf{X}_t^0. \tag{2.8'}$$

$a = 1, b = 5, c = 1, \sigma = \varkappa = -2, m = n = 2, x_0 = 0, T = 5, R = T/0.01.$

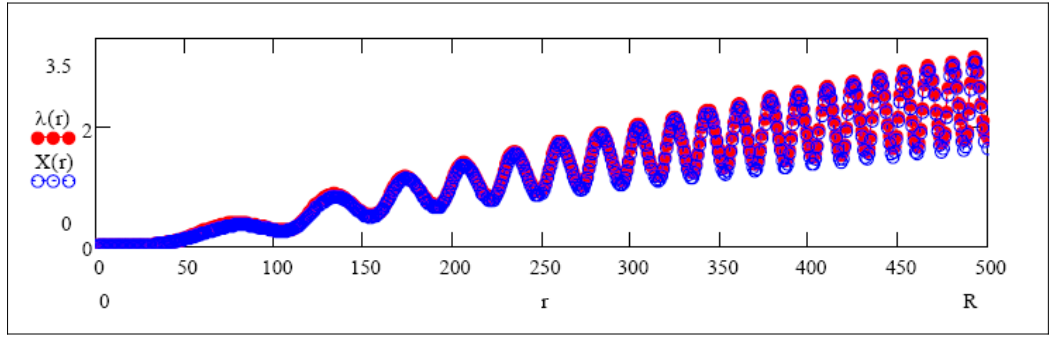

The solution Eq.(8) in comparison with a coresponding solution X(t) ODE (8′)

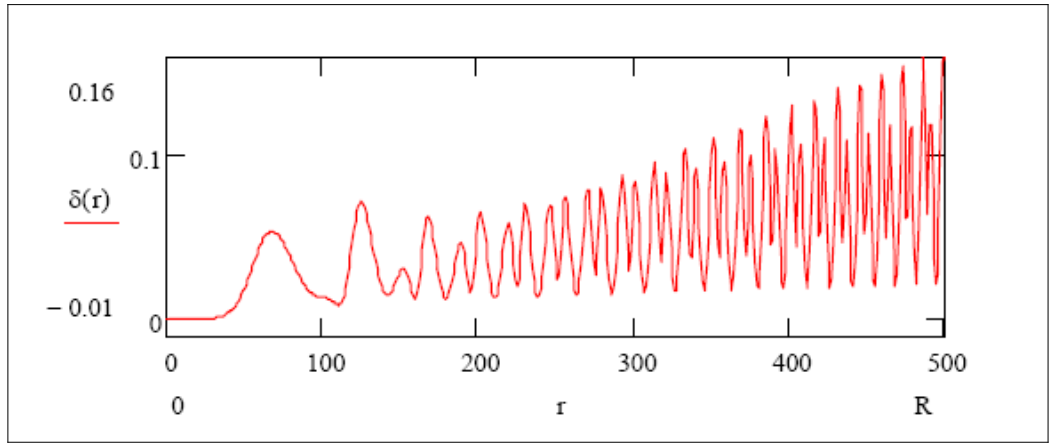

Here $\delta(t) \triangleq \liminf_{\varepsilon \to 0} \mathbf{E}\left[\|\mathbf{X}_t^\varepsilon - \mathbf{X}_t^0\|^2\right].$ Let us consider now an $m$-persons Colombeau stochastic differential game $\mathbf{CDG}_{m;T}(\mathbf{f}, \mathbf{g}, \mathbf{y})$ with nonlinear dynamics

$$(\dot{\mathbf{X}}_{\epsilon,\varepsilon}(t,\omega))_\epsilon = \mathbf{f}_\epsilon(t, (\mathbf{X}_{\epsilon,\varepsilon}(t,\omega))_\epsilon, \boldsymbol{\alpha}(t)) + \sqrt{\varepsilon}\, d\mathbf{W}(t); \forall t : (x_\epsilon(t))_\epsilon \in \widetilde{\mathbb{R}}^n,$$

$$x(0) = x_0, t \in [0,T], \epsilon \in (0,1], \varepsilon \ll 1, \quad (2.9)$$

$$\boldsymbol{\alpha}(t) = (\alpha_1(t), \ldots, \alpha_m(t)), \alpha_i(t) \in U_i \subsetneq \mathbb{R}^{k_i}, i = 1, \ldots, m.$$

Here $t \to \alpha_i(t)$ is the control chosen by the $i$-th player, within a set of admissible control values $U_i$, and the payoff of the $i$-th player is

$$(\bar{\mathbf{J}}_{\epsilon,i})_\epsilon = \mathbf{E}\left[\int_0^T (g_{\epsilon,i}(\mathbf{X}_{\epsilon,\varepsilon}(t,\omega);\boldsymbol{\alpha}(t)))_\epsilon dt\right] + \mathbf{E}\left[\sum_{i=1}^n [(X_{\epsilon,\varepsilon;i}(T,\omega))_\epsilon - y_i]^2\right] \quad (2.10)$$

where $\mathbf{y} = (y_1, \ldots, y_n)$ and $t \mapsto x(t,\omega)$ is the trajectory of the Eq.(2.9).

**Theorem.2.** For any solution $\{(\mathbf{X}_t^{\epsilon,\varepsilon})_\epsilon, \check{\alpha}(t)\} = \{((X_{1,t}^{\epsilon,\varepsilon})_\epsilon, \ldots, (X_{n,t}^{\epsilon,\varepsilon})_\epsilon), (\check{\alpha}_1(t), \ldots, \check{\alpha}_m(t))\}$ of the dissipative $\mathbf{CDG}_{m;T}(\mathbf{f}, \mathbf{0}, \mathbf{y})$ and $\mathbb{R}$ valued parameters $\lambda_1, \ldots, \lambda_n$, there exists constant $C' \in \mathbb{R}_+$ such that:

$$\liminf_{\varepsilon \to 0} \mathbf{E}\left[\|\mathbf{X}_t^\varepsilon - \lambda\|^2\right] \leq C' \|\mathbf{U}(t,\lambda)\|^2. \quad (2.11)$$

Where $\mathbf{U}(t,\lambda) = (u_1(t,\lambda), \ldots, u_n(t,\lambda))$ is a trajectory of the *corresponding linear master game:*

$$\dot{\mathbf{U}}(t,\lambda) = \mathbf{J}[\mathbf{f}(\lambda, \check{\alpha}(t,\lambda))]\mathbf{U} + \mathbf{f}(\lambda, \check{\alpha}(t,\lambda)), \mathbf{U}(0,\lambda) = x_0 - \lambda,$$

$$\bar{\mathbf{J}}_i = \|\mathbf{U}(T)\|^2, \quad (2.12)$$

$$\left(\min_{\alpha_i(t)} \left(\max_{\alpha_j(t), j \neq i} \bar{\mathbf{J}}_i\right)\right).$$

**Theorem.3.** Suppose that: **(i)** there exists function $\lambda : [0,T] \to \mathbb{R}^n$ such that $\forall t (t \in [0,T])$ the next condition is satisfied

$$\mathbf{U}(t, \lambda(t)) = 0. \quad (2.13)$$

**(ii)** feedback optimal control of the $i$-th player for a linear master game (2.12) with any fixed parameter $\lambda \in \mathbb{R}^n$ is

$$\check{\alpha}_i(t,\lambda) = \check{\alpha}_i[t, u_1(t), \ldots, u_n(t); \lambda_1, \ldots, \lambda_n], \qquad (2.14)$$

where $i = 1, \ldots, m$.

Then feedback optimal control of the $i$-th player for the $m$-persons nonlinear game (2.9)-(2.10) is

$$\check{\alpha}_i(t,\lambda) = \check{\alpha}_i[t, x_1(t), \ldots, x_n(t); x_1(t), \ldots, x_n(t)]. \qquad (2.15)$$

**Example.2. Nonlinear $2$-persons differential game.**

$$\dot{x}_1 = x_2, \dot{x}_2 = -\kappa x_2^3 + \alpha_1(t) + \alpha_2(t), \kappa > 0, \ t \in [0,T],$$

$$x_1(0) = x_{01}, x_2(0) = x_{02};$$

$$\alpha_1(t) \in [-\rho_1, \rho_1], \alpha_2(t) \in [-\rho_2, \rho_2], \qquad (2.16)$$

$$\mathbf{J}_i(T) = x_1^2(T) + x_2^2(T), i = 1, 2.$$

Optimal control problem of the first player is

$$\min_{\alpha_1(t) \in [-\rho_1, \rho_2]} \left( \max_{\alpha_2(t) \in [-\rho_2, \rho_2]} \mathbf{J}_1(T) \right) =$$

$$\min_{\alpha_1(t) \in [-\rho_1, \rho_2]} \left( \max_{\alpha_2(t) \in [-\rho_2, \rho_2]} [x_1^2(T) + x_2^2(T)] \right) \qquad (2.17)$$

and optimal control problem of the second player is

$$\max_{\alpha_2(t)\in[-\rho_u,\rho_u]} \left( \min_{\alpha_1(t)\in[-\rho_1,\rho_1]} \mathbf{J}_2(T) \right) =$$

$$\max_{\alpha_2(t)\in[-\rho_u,\rho_u]} \left( \min_{\alpha_1(t)\in[-\rho_1,\rho_1]} [x_1^2(T) + x_2^2(T)] \right). \tag{2.18}$$

From nonlinear 2-persons differential game given by Eqs.(2.16)-(2.18) we obtain corresponding linear master game:

$$\dot{u}_1 = u_2 + \lambda_2,$$

$$\dot{u}_2 = -3\kappa\lambda_2^2 u_2 - \kappa\lambda_2^3 + \check{\alpha}_1(t) + \check{\alpha}_2(t),$$

$$u_1(0) = x_{01} - \lambda_1, u_2(0) = x_{02} - \lambda_2, \tag{2.19}$$

$$\check{\alpha}_1(t) \in [-\rho_1,\rho_1], \check{\alpha}_2(t) \in [-\rho_2,\rho_2],$$

$$\bar{\mathbf{J}}_i = u_1^2(T) + u_2^2(T), i = 1,2;$$

optimal control problem of the first player is

$$\min_{\check{\alpha}_1(t)\in[-\rho_1,\rho_2]} \left( \max_{\check{\alpha}_2(t)\in[-\rho_2,\rho_2]} [u_1^2(T) + u_2^2(T)] \right) \tag{2.20}$$

and optimal control problem of the second player:

$$\max_{\check{\alpha}_2(t)\in[-\rho_u,\rho_u]} \left( \min_{\check{\alpha}_1(t)\in[-\rho_1,\rho_1]} [u_1^2(T) + u_2^2(T)] \right). \tag{2.21}$$

Having solved by standard way [7]-[8] linear master game (2.19) one obtain optimal feedback control of the first player:

$$\check{\alpha}_1(t) = \check{\alpha}_1[t, u_1(t), u_2(t)] =$$

$$-\rho_1 \mathbf{sign}\{u_1(t) + [(T-t) + \exp[-3\kappa u_2^2(t)(T-t)]]u_2(t)\} \qquad (2.22)$$

and optimal feedback control of the second player:

$$\check{\alpha}_2(t) = \check{\alpha}_2[t, u_1(t), u_2(t)] =$$

$$\rho_2 \mathbf{sign}\{u_1(t) + [(T-t) + \exp[-3\kappa u_2^2(t)(T-t)]]u_2(t)\}. \qquad (2.23)$$

Using **Theorem 3** we obtain complete solution of the nonlinear 2-persons differential game given by Eqs.(14-16) in the next form. Optimal feedback control of the first player is

$$\alpha_1(t) = \check{\alpha}_1[t, x_1(t), x_2(t)] =$$

$$-\rho_1 \mathbf{sign}\{x_1(t) + [(T-t) + \exp[-3\kappa x_2^2(t)(T-t)]]x_2(t)\} \qquad (2.24)$$

and optimal feedback control of the second player is

$$\alpha_2(t) = \check{\alpha}_2[t, x_1(t), x_2(t)] =$$

$$\rho_2 \mathbf{sign}\{x_1(t) + [(T-t) + \exp[-3\kappa x_2^2(t)(T-t)]]x_2(t)\}. \qquad (2.25)$$

For numerical simulation from Eq.(2.16) and Eqs.(2.24-(2.25) we obtain ODE:

$$\dot{x}_1(t) = x_2(t),$$

$$\dot{x}_2(t) = -\kappa x_2^3(t) - \quad (2.26)$$

$$\rho_1 \cdot \text{sign}[x_1(t) + [(T-t) + \exp[-3\kappa x_2^2(t)(T-t)]]x_2(t)] + \alpha_2(t).$$

**Numerical simulation:** $\kappa = 1, \rho_1 = 400, A = 100, \omega = 5, x_1(0) = 300m, x_2(0) = 30m/\text{sec}, T = 80\text{sec}, \alpha_2(t) = A\sin^2(\omega \cdot t).$

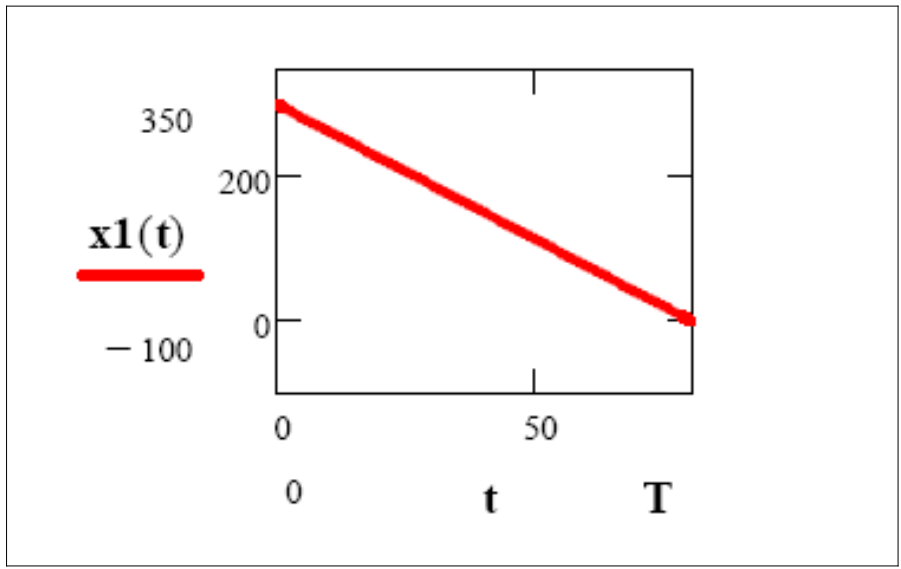

Optimal trajectory : $x_1(t). x_1(T) = 0.4m$

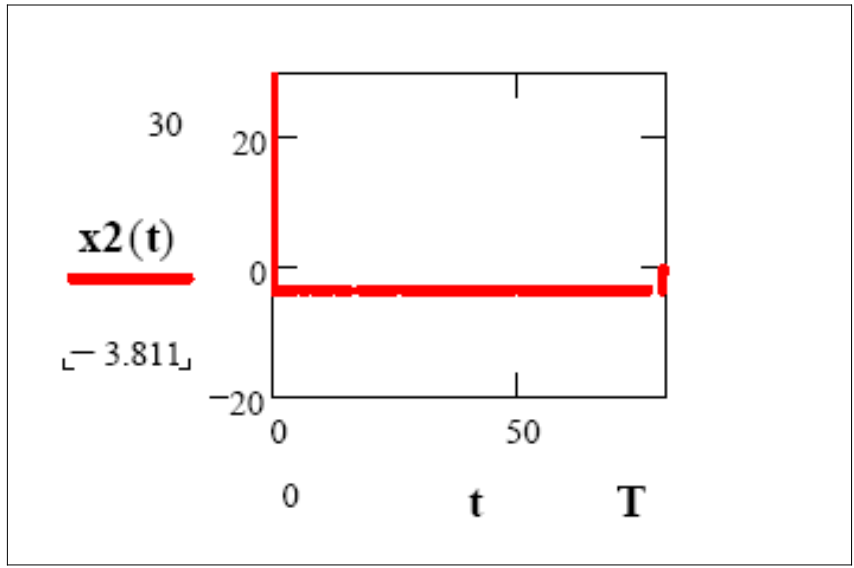

Optimal velocity : $x_2(t). x_2(T) = -0.4 m/\sec$

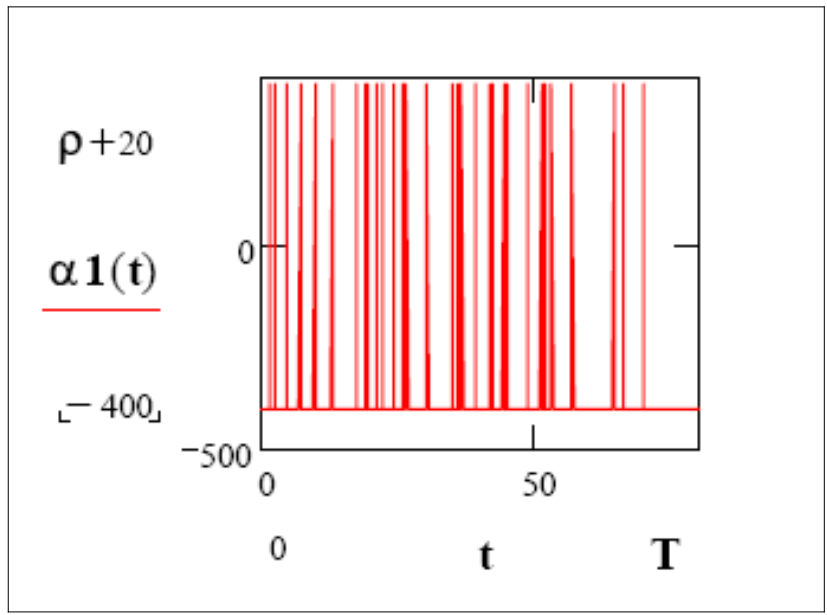

Optimal control of the first player

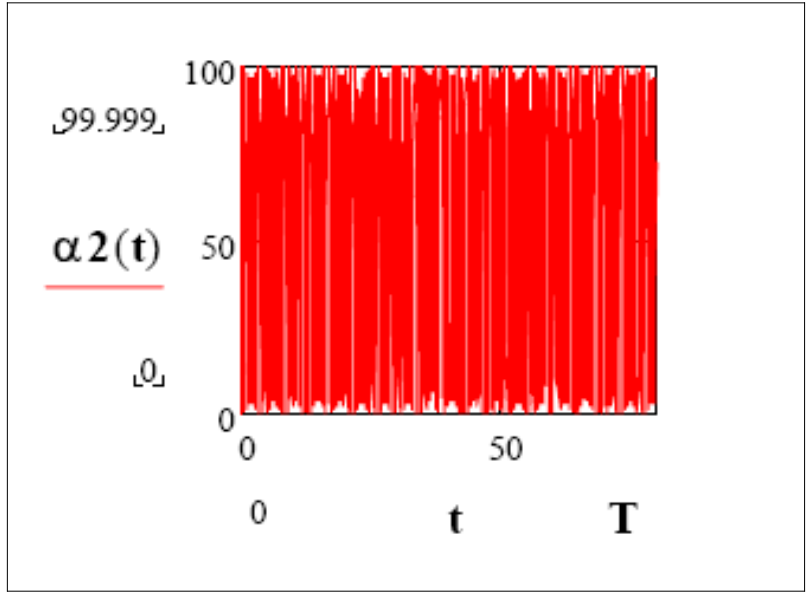

Control of the second player.

Let us consider now an $m$-persons Colombeau stochastic differential game $\mathbf{CDG}_{m;T}(\mathbf{f},\mathbf{g},\mathbf{y})$ with imperfect information and nonlinear dynamics

$$(\dot{\mathbf{X}}_{\epsilon,\varepsilon}(t,\omega))_{\epsilon} = \mathbf{f}_{\epsilon}(t,(\mathbf{X}_{\epsilon,\varepsilon}(t,\omega))_{\epsilon},\boldsymbol{\alpha}(t);\boldsymbol{\beta}(t)) + \sqrt{\varepsilon}\,d\mathbf{W}(t); \forall t : (x_{\epsilon}(t))_{\epsilon} \in \widetilde{\mathbb{R}}^{n},$$

$$x(0) = x_0, t \in [0,T], \epsilon \in (0,1], \varepsilon \ll 1, \qquad (2.27)$$

$$\boldsymbol{\alpha}(t) = (\alpha_1(t;\boldsymbol{\beta}(t)),\ldots,\alpha_m(t;\boldsymbol{\beta}(t))), \alpha_i(t;\boldsymbol{\beta}(t)) \in U_i \subsetneq \mathbb{R}^{k_i}, i = 1,\ldots,m.$$

Here $\boldsymbol{\beta}(t)$ is an uncertainty, $t \to \alpha_i(t;\boldsymbol{\delta}(t))$ is the control chosen by the $i$-th player, within a set of admissible control values $U_i$, and the payoff of the $i$-th player is

$$(\bar{\mathbf{J}}_{\epsilon,i})_{\epsilon} = \mathbf{E}\left[\int_0^T (g_{\epsilon,i}(\mathbf{X}_{\epsilon,\varepsilon}(t,\omega);\boldsymbol{\alpha}(t;\boldsymbol{\delta}(t))))_{\epsilon}dt\right] + \mathbf{E}\left[\sum_{i=1}^n [(X_{\epsilon,\varepsilon;i}(T,\omega))_{\epsilon} - y_i]^2\right] \qquad (2.28)$$

where $\mathbf{y} = (y_1,\ldots,y_n)$ and $t \mapsto x(t,\omega)$ is the trajectory of the Eq.(2.27).

**Theorem.4.** For any solution $\{(\mathbf{X}_t^{\epsilon,\varepsilon})_\epsilon, \check{\alpha}(t)\} = \{\left((X_{1,t}^{\epsilon,\varepsilon})_\epsilon, \ldots, (X_{n,t}^{\epsilon,\varepsilon})_\epsilon\right), (\check{\alpha}_1(t), \ldots, \check{\alpha}_m(t))\}$ dissipative $\mathbf{CIDG}_{m;T}(\mathbf{f}, \mathbf{0}, \mathbf{y}, \boldsymbol{\beta})$ and $\mathbb{R}$ valued parameters $\lambda_1, \ldots, \lambda_n$, there exists constant $\mathbf{C}' = C' \in \mathbb{R}_+$, such that:

$$\liminf_{\varepsilon \to 0} \mathbf{E}\left[ \|\mathbf{X}_t^\varepsilon - \lambda\|^2 \right] \le \mathbf{C}' \|\mathbf{U}(t, \lambda; \boldsymbol{\beta}(t))\|^2. \tag{2.29}$$

Where $\mathbf{U}(t, \lambda; \boldsymbol{\beta}(t)) = (U_1(t, \lambda; \boldsymbol{\beta}(t)), \ldots, U_n(t, \lambda; \boldsymbol{\beta}(t)))$ the trajectory of the *corresponding linear master game with imperfect information*

$$\dot{\mathbf{U}}(t, \lambda; \boldsymbol{\beta}(t)) =$$

$$\mathbf{J}[\mathbf{f}(\lambda, \check{\alpha}(t, \lambda; \boldsymbol{\beta}(t)))]\mathbf{U}(\mathbf{t}, \lambda; \boldsymbol{\beta}(\mathbf{t})) + \mathbf{f}(\lambda, \check{\alpha}(t, \lambda; \boldsymbol{\beta}(t))), \mathbf{U}(0, \lambda) = x_0 - \lambda,$$

$$\bar{\mathbf{J}}_i = \|\mathbf{U}(T)\|^2, \tag{2.30}$$

$$\left(\min_{\alpha_i(t)}\left(\max_{\alpha_j(t), j \ne i} \bar{\mathbf{J}}_i\right)\right).$$

**Theorem.5.** Suppose that: **(i)** there exists function $\lambda : [0, T] \to \mathbb{R}^n$ such that $\forall t (t \in [0, T])$ the next condition is satisfied

$$\mathbf{U}(t, \lambda(t), \boldsymbol{\beta}(t)) = 0. \tag{2.31}$$

**(ii)** feedback optimal control of the $i$-th player for a linear master game (2.30) with any fixed parameter $\lambda \in \mathbb{R}^n$ is

$$\check{\alpha}_i(t, \lambda) = \check{\alpha}_i[t, u_1(t), \ldots, u_n(t); \lambda_1, \ldots, \lambda_n; \boldsymbol{\beta}(t)], \tag{2.32}$$

where $i = 1,\ldots,m$.

Then feedback optimal control of the $i$-th player for the $m$-persons nonlinear game (2.27)-(2.28) is

$$\check{\alpha}_i(t,\lambda) = \check{\alpha}_i[t,x_1(t),\ldots,x_n(t);x_1(t),\ldots,x_n(t);\delta(t),\beta(t)]. \qquad (2.33)$$

**Example.3.** Nonlinear 2-persons differential game with imperfect measurements.

$$\dot{x}_1(t) = x_2(t),$$

$$\dot{x}_2(t) = -\kappa_1 x_2^3(t) + \kappa_2 x_2^2(t) + \alpha_1[t,x_1(t),x_2(t)+\beta(t)] + \alpha_2[t,x_1(t),x_2(t)],$$

$$\kappa_1 > 0, \kappa_2 \in (-\infty,+\infty), \alpha_1(t) \in [-\rho_1,\rho_1], \alpha_2(t) \in [-\rho_2,\rho_2], \qquad (2.34)$$

$$\bar{\mathbf{J}}_i = x_1^2(T) + \dot{x}_1^2(T), i = 1,2,$$

$$\left(\min_{\alpha_i(t)} \left(\max_{\alpha_j(t), j\neq i} \bar{\mathbf{J}}_i\right)\right).$$

From Eq.(28) one obtain corresponding linear master game:

$$\dot{u}_1 = u_2 + \lambda_2,$$

$$\dot{u}_2 = -(3\kappa_1\lambda_2^2 - 2\kappa_2\lambda_2)u_2 - \kappa_1\lambda_2^3 +$$

$$+\kappa_2\lambda_2^2 + \check{\alpha}_1[t, u_1(t), u_2(t) + \beta(t)] + \check{\alpha}_2[t, u_1(t), u_2(t)], \qquad (2.35)$$

$$\check{\alpha}_1(t) \in [-\rho_1, \rho_1], \check{\alpha}_2(t) \in [-\rho_2, \rho_2],$$

$$\bar{J}_i = u_1^2(T) + u_2^2(T), i = 1, 2.$$

Having solved by standard way linear master game given by Eq.(32) one obtain local optimal feedback control of the first player [2] chapt.IV.3:

$$\check{\alpha}_{1,\mathbf{loc}}(t) = -\rho_1 \mathbf{sign}[u_1(t) + (t_{n+1} - t)(u_2(t) + \beta(t))] \qquad (2.36)$$

and local optimal feedback control of the second player:

$$\check{\alpha}_{2,\mathbf{loc}}(t) = \rho_2 \mathbf{sign}[u_1(t) + (t_{n+1} - t)u_2(t)]. \qquad (2.37)$$

Finally we obtain global optimal control in the next form [2] chapt.IV.3:

$$\check{\alpha}_1(t) = -\rho_1 \mathbf{sign}[u_1(t) + \Theta_\tau(t)(u_2(t) + \beta(t))],$$

$$\check{\alpha}_2(t) = \rho_2 \mathbf{sign}[u_1(t) + \Theta_\tau(t)u_2(t)]. \qquad (2.38)$$

Here

$$\Theta_\tau(t) = \theta_\tau(\eta_\tau(t)), \theta_\tau(t) \triangleq \tau - t,$$

$$\eta_\tau(t) \triangleq t - \left(\mathbf{ceil}\left(\frac{t}{\tau}\right) - 1\right),$$

(2.39)

where $\mathbf{ceil}(x)$ is a part-whole number of a number $x \in \mathbb{R}$.

$$\alpha_1(t) = -\rho_1 \mathbf{sign}[x_1(t) + \Theta_\tau(t)(x_2(t) + \beta(t))],$$

$$\alpha_2(t) = \rho_2 \mathbf{sign}[x_1(t) + \Theta_\tau(t)x_2(t)].$$

(2.40)

Thus for numerical simulation we obtain ODE:

$$\dot{x}_1 = x_2,$$
$$\dot{x}_2 = -\kappa_1 x_2^3 + \kappa_2 x_2^2 - \rho_1 \cdot \mathbf{sign}[x_1(t) + \Theta_\tau(t)(x_2(t) + \beta(t))] +$$
$$\rho_2 \mathbf{sign}[x_1(t) + \Theta_\tau(t)x_2(t)].$$

(2.41)

**Numerical simulation.** Game with imperfect measurements: red curves $\mathbf{x1(t), x2(t)}$. Classical game: blue curves $\mathbf{y1(t), y2(t)}$. $\beta(t) = A\sin^2(\omega \cdot t)$

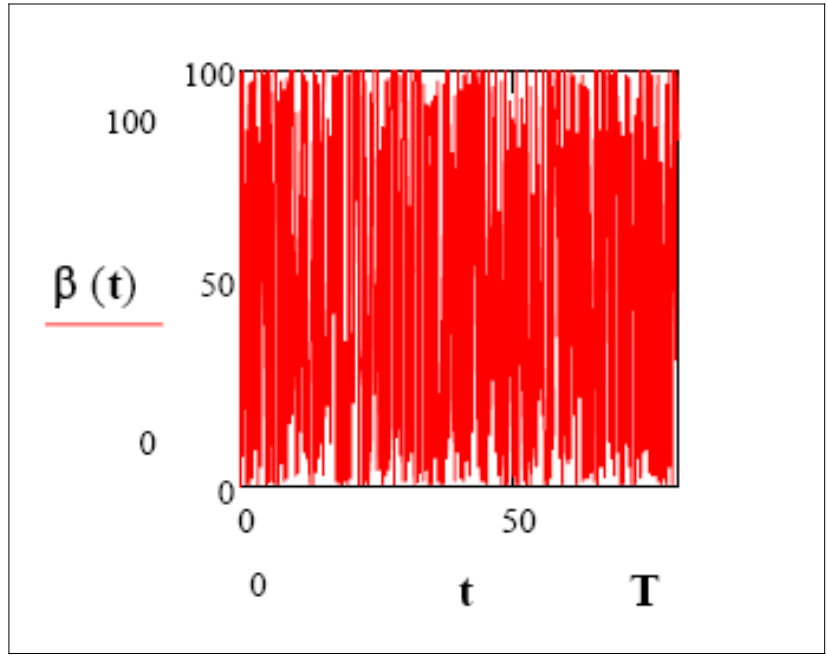

Uncertainty of speed measurements $\beta(t)$.

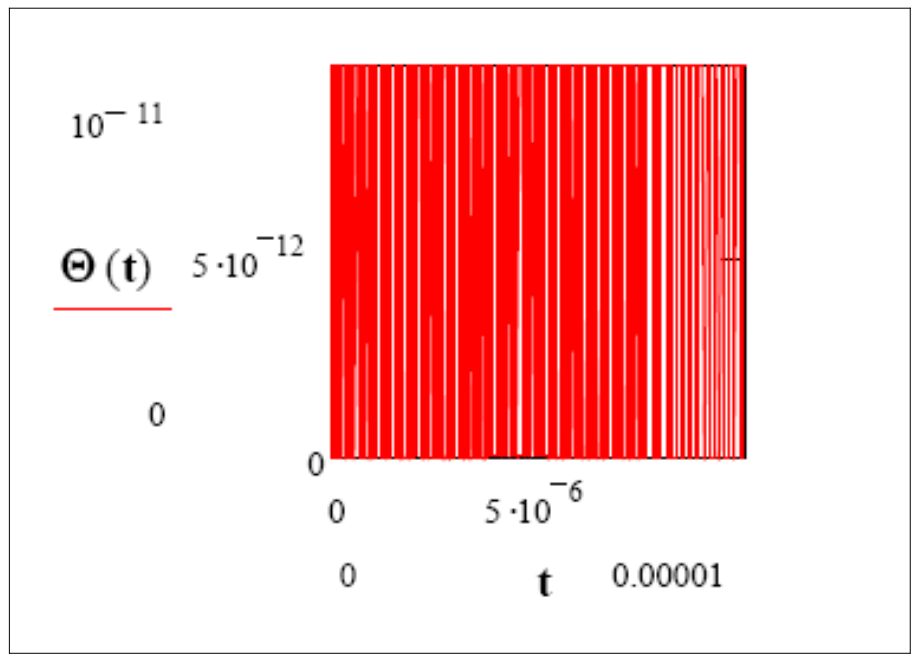

Cutting function $\Theta_\tau(t). \tau = 10^{-11}$.

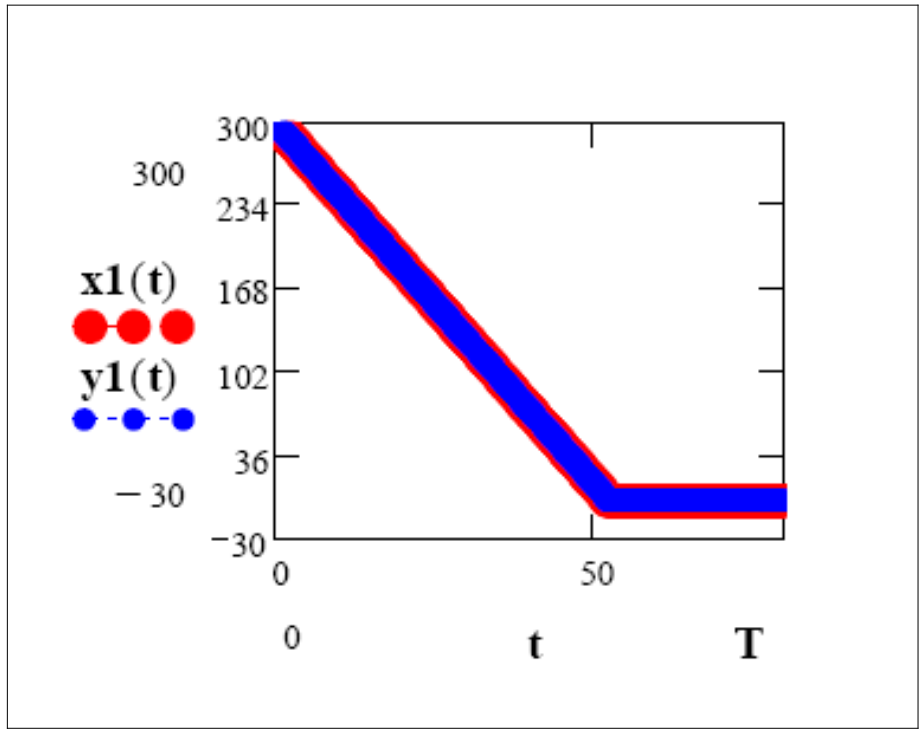

Optimal trajectory

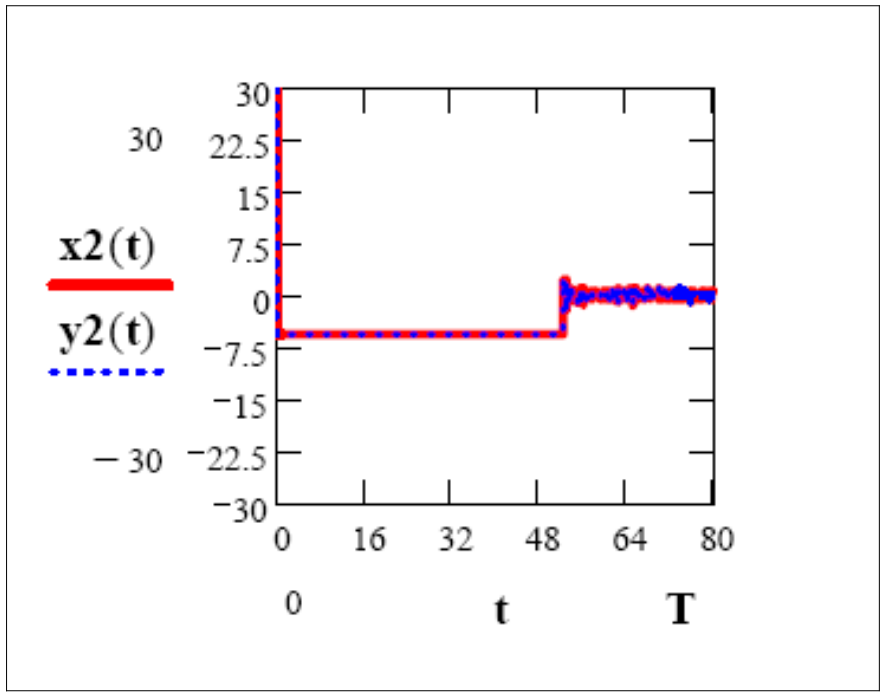

Optimal velocity.

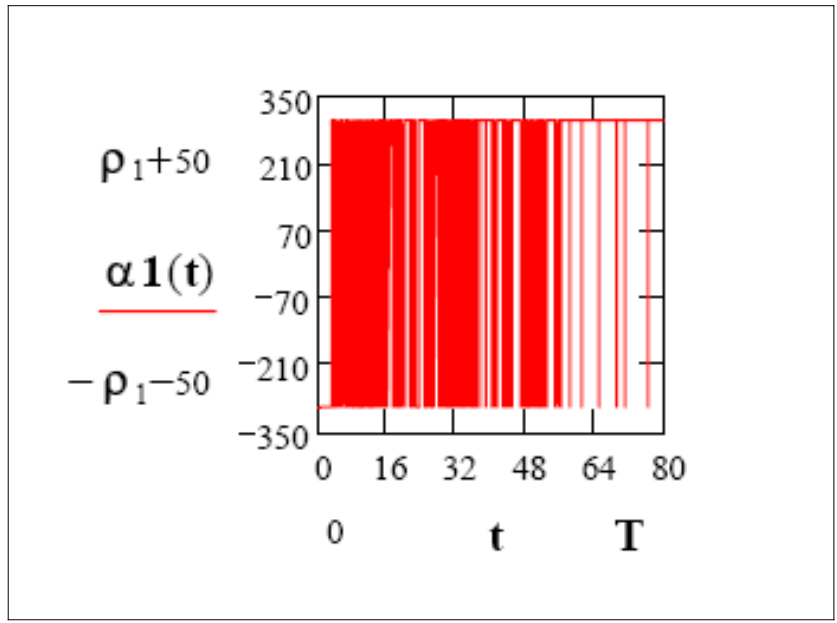

Optimal control of the first player

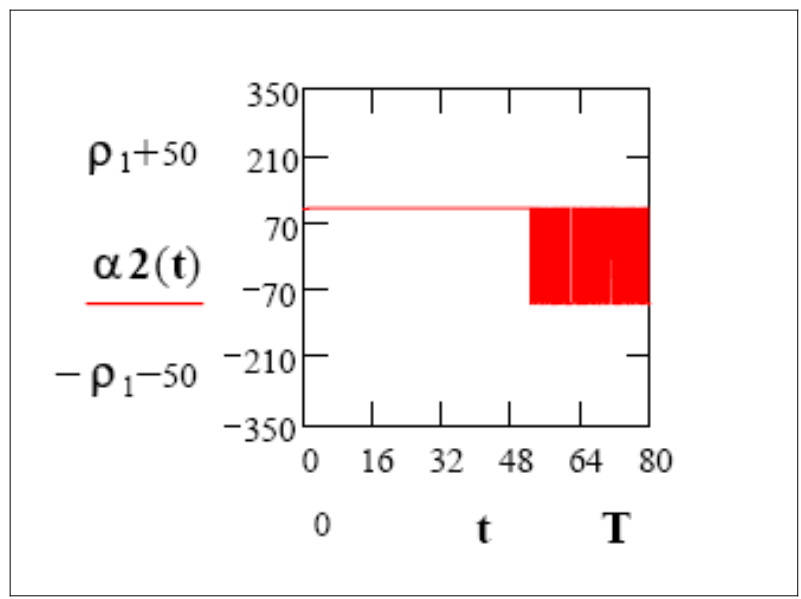

Optimal control of the second player

**3. Supporting Technical Analysis.** Let us consider optimal control problem from Example.1. Corresponding Bellman type equation is

$$\min_{\alpha_1\in[-\rho_1,\rho_1]} \left( \max_{\alpha_2\in[-\rho_2,\rho_2]} \left[ \frac{\partial V}{\partial t} + \frac{\partial V}{\partial x_1}x_2 + \frac{\partial V}{\partial x_2}(-x_2^3 + \alpha_1 + \alpha_2) \right] \right) = 0, V(T,x)$$

$$= (x_1^2 + x_2^2).$$
(3.39)

Complete constructing the exact analytical solution for PDE (39) is a complicated unresolved classical problem,because PDE (39) is not amenable to analytical treatments.Even the theorem of existence classical solution for boundary Problems such (39) is not proved. Thus even for simple cases a problem of construction feedback optimal control by the associated Bellman equation complicated numerical technology or principal simplification is needed [4].Hovewer as one can see complete constructing feedback optimal control from Theorem (3) is simple.

### 4.Homing missile guidance with imperfect measurements capable to defeat in conditions of hostile active radio-electronic jamming.

Homing missile guidance strategies (guidance laws) dictate the manner in which the missile will guide to intercept, or rendezvous with, the target. The feedback nature of homing guidance allows the guided missile (or, more generally, the pursuer) to tolerate some level of (sensor) measurement uncertainties, errors in the assumptions used to model the engagement (e.g., unanticipated target maneuver), and errors in modeling missile capability (e.g., deviation of actual missile speed of response to guidance commands from the design assumptions).Nevertheless, the selection of a guidance strategy and its subsequent mechanization are crucial design factors that can have substantial impact on guided missile performance. Key drivers to guidance law design include the type of targeting sensor to be used (passive IR, active or semi-active RF, etc.), accuracy of the targeting and inertial measurement unit (IMU) sensors, missile maneuverability, and, finally yet important, the types of targets to be engaged and their associated maneuverability levels.

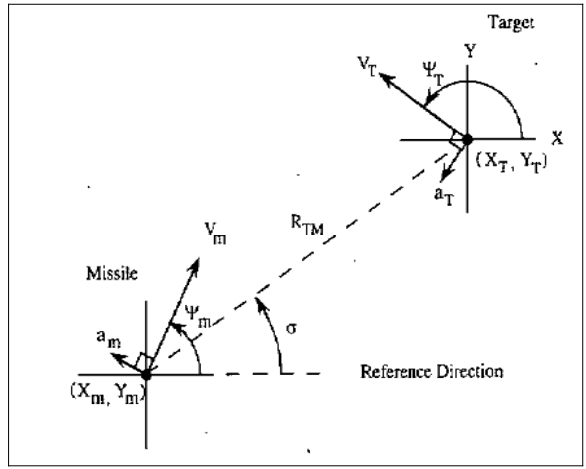

Fig4.1.Planaimtercegeometry

## I.Homing missile guidance with a perfect measurements

Figure 4.1 shows the ntercept geometry of a missile in planar pursuit of a target. Taking the origin of the reference frame to be the instantaneous position of the
missile, the equation of motion in polar form are:

$$\ddot{R} = R\dot{\sigma}^2 + a_m^r(t) + a_T^r(t),$$

$$R\ddot{\sigma} + 2\dot{R}\dot{\sigma} = a_M^\tau(t) + a_T^\tau(t),$$

$$a_m^r(t) = a_m(t)\sin(\Psi_m - \sigma), a_m^r(t) = a_m(t)\sin(\Psi_m - \sigma),$$

$$a_m^\tau(t) = a_m(t)\cos(\Psi_T - \sigma), a_T^\tau(t) = a_T(t)\cos(\Psi_T - \sigma), \qquad (4.1)$$

$$\dot{\Psi}_m = \frac{a_m}{V_m}, \dot{\Psi}_T = \frac{a_T}{V_T},$$

$$a_M^r(t) \in [-\bar{a}_m^r, \bar{a}_m^r], a_T^r(t) \in [-\bar{a}_T^r, \bar{a}_T^r],$$

$$a_M^\tau(t) \in [-\bar{a}_m^\tau, \bar{a}_m^\tau], a_T^\tau(t) \in [-\bar{a}_T^\tau, \bar{a}_T^\tau]$$

1.The variable $R = R(t)$ denotes the target-to-missile range $R_{TM}(t)$.
2.The variable $\sigma = \sigma(t)$ denotes the line-of-sight angle (LOS) i.e.,the angle between the constant reference direction and target-to-missile direction.

4. The variable $a_M^\tau(t)$ denotes the missiles tangent acceleration,i.e.missile acceleration along direction which perpendicularly to line-of-sight direction.

5. The variable $a_M^r(t)$ denotes the missile acceleration along target-to-missile direction.

6. The variable $a_T^\tau(t)$ denotes the target tangent acceleration.

7. The variable $a_T^r(t)$ denotes the target acceleration along target-to-missile direction.

Using the replacement $\dot{z} = R\dot{\sigma}$ from Eq.(4.1) one obtain:

$$\dot{\sigma} = \frac{\dot{z}}{R},$$

$$\ddot{z} = \dot{R}\dot{\sigma} + R\ddot{\sigma}, R\ddot{\sigma} = \ddot{z} - \dot{R}\dot{\sigma}, \tag{4.2}$$

$$R\ddot{\sigma} + 2\dot{R}\dot{\sigma} = \ddot{z} + \dot{R}\dot{\sigma} = \ddot{z} + \frac{\dot{R}\dot{z}}{R}.$$

Substitution Eq.(4.2) into Eq.(4.1) gives:

$$\ddot{R} = \frac{\dot{z}^2}{R} + a_M^r(t) + a_T^r(t),$$

$$a_M^r(t) \in [-\bar{a}_M^r, \bar{a}_M^r], a_T^r(t) \in [-\bar{a}_T^r, \bar{a}_T^r].$$

$$\ddot{z} = -\frac{\dot{R}\dot{z}}{R} + a_M^\tau(t) + a_T^\tau(t), \tag{4.3}$$

$$a_M^\tau(t) \in [-\bar{a}_M^\tau, \bar{a}_M^\tau], a_T^\tau(t) \in [-\bar{a}_T^\tau, \bar{a}_T^\tau].$$

Using the regularization of the Eq.(4.1) we obtain:

$$\ddot{R} = \frac{\dot{z}^2}{R+\varepsilon} + a_M^r(t) + a_T^r(t),$$

$$a_M^r(t) \in [-\bar{a}_M^r, \bar{a}_M^r], a_T^r(t) \in [-\bar{a}_T^r, \bar{a}_T^r].$$

(4.4)

$$\ddot{z} = -\frac{\dot{R}\dot{z}}{R+\varepsilon} + a_M^\tau(t) + a_T^\tau(t),$$

$$a_M^\tau(t) \in [-\bar{a}_M^\tau, \bar{a}_M^\tau], a_T^\tau(t) \in [-\bar{a}_T^\tau, \bar{a}_T^\tau].$$

Here $\varepsilon \in (0,1]$. Using the replacement $\dot{R} = V_m, \dot{z} = w$ into Eq.(4.2.24) one obtain:

$$\dot{R} = V_m,$$

$$\dot{V}_m = \frac{w^2}{R+\varepsilon} + a_M^r(t) + a_T^r(t),$$

$$a_M^r(t) \in [-\bar{a}_M^r, \bar{a}_M^r], a_T^r(t) \in [-\bar{a}_T^r, \bar{a}_T^r];$$

(4.5)

$$\dot{z} = w,$$

$$\dot{w} = -\frac{V_r w}{R+\varepsilon} + a_M^\tau(t) + a_T^\tau(t),$$

$$a_M^\tau(t) \in [-\bar{a}_M^\tau, \bar{a}_M^\tau], a_T^\tau(t) \in [-\bar{a}_T^\tau, \bar{a}_T^\tau].$$

Let us consider the optimal control problem:

$$\dot{R} = V_m,$$

$$\dot{V}_m = \frac{w^2}{R+\varepsilon} + a_M^r(t) + a_T^r(t) - \kappa_1 V_m^3,$$

$$\check{a}_M^r(t) \in [-\bar{a}_M^r, \bar{a}_M^r], a_T^r(t) \in [-\bar{a}_T^r, \bar{a}_T^r], 0 < \kappa_1,$$

$$\dot{z} = w, \qquad (4.6)$$

$$\dot{w} = -\frac{V_r w}{R+\varepsilon} + a_m^\tau(t) + a_T^\tau(t) - \kappa_2 w^3, 0 < \kappa_2,$$

$$a_m^\tau(t) \in [-\bar{a}_m^\tau, \bar{a}_m^\tau], a_T^\tau(t) \in [-\bar{a}_T^\tau, \bar{a}_T^\tau],$$

$$\mathbf{J}_i = R^2(t_1) + z^2(t_1), i = 1, 2.$$

Optimal control problem of the first player is:

$$\mathbf{J}_1 = \min_{a_M^r(t) \in [-\bar{a}_M^r, \bar{a}_M^r], a_T^r(t) \in [-\bar{a}_T^r, \bar{a}_T^r],} \left\{ \max_{a_T^r(t) \in [-\bar{a}_T^r, \bar{a}_T^r], a_T^\tau(t) \in [-\bar{a}_T^\tau, \bar{a}_T^\tau]} [R^2(t_1) + z^2(t_1)] \right\}. \qquad (4.7)$$

Optimal control problem of the second player is:

$$\mathbf{J}_2 = \max_{a_T^r(t) \in [-\bar{a}_T^r, \bar{a}_T^r], a_T^\tau(t) \in [-\bar{a}_T^\tau, \bar{a}_T^\tau]} \left\{ \max_{a_M^r(t) \in [-\bar{a}_M^r, \bar{a}_M^r], a_T^\tau(t) \in [-\bar{a}_T^r, \bar{a}_T^r].} [R^2(t_1) + z^2(t_1)] \right\}. \qquad (4.8)$$

We assume now for shortness only that (1) $\dot{z}^2/R = R\dot{\sigma}^2 \simeq 0$, (2) $V_m \dot{z}/R = V_m \dot{\sigma} \simeq 0$. Let us consider now the next optimal control problem:

$$\dot{r} = v_m,$$

$$\dot{v}_m = a_m^r(t) + a_T^r(t) - \kappa_1 V_m^3,$$

$$a_m^r(t) \in [-\bar{a}_m^r, \bar{a}_m^r], a_T^r(t) \in [-\bar{a}_T^r, \bar{a}_T^r], 0 < \kappa_1.$$

(4.9)

$$\dot{w} = a_m^\tau(t) + a_T^\tau(t) - \kappa_2 w^3,$$

$$a_m^\tau(t) \in [-\bar{a}_m^\tau, \bar{a}_m^\tau], a_T^\tau(t) \in [-\bar{a}_T^\tau, \bar{a}_T^\tau], 0 < \kappa_2.$$

$$\mathbf{J}_i = R^2(t_1) + z^2(t_1), i = 1, 2.$$

Using solution which given in [2] chapt.IV.2 Eq.(4.2.13) we obtain global quasy optimal control in the next form:

$$a_m^r(t) = -\bar{a}_m^r \mathbf{sign}[R(t) + \Theta_\tau(t)(V_m(t))],$$

(4.10)

$$a_m^\tau(t) = -\bar{a}_m^\tau \mathbf{sign}[z(t) + \Theta_\tau(t)w(t)].$$

Thus for numerical simulation we obtain ODE:

$$\dot{R} = V_m,$$

$$\dot{V}_m = \frac{w^2}{R+\varepsilon} - \bar{a}_m^r \operatorname{sign}[R(t) + \Theta_\tau(t)(V_m(t))] + a_T^r(t) - \kappa_1 V_m^3,$$

$$\check{a}_M^r(t) \in [-\bar{a}_M^r, \bar{a}_M^r], a_T^r(t) \in [-\bar{a}_T^r, \bar{a}_T^r],$$

$$\dot{z} = w, \tag{4.11}$$

$$\dot{w} = -\frac{V_r w}{R+\varepsilon} - \bar{a}_m^\tau \operatorname{sign}[z(t) + \Theta_\tau(t)w(t)] + a_T^\tau(t) - \kappa_2 w^3,$$

$$a_m^\tau(t) \in [-\bar{a}_m^\tau, \bar{a}_m^\tau], a_T^\tau(t) \in [-\bar{a}_T^\tau, \bar{a}_T^\tau],$$

$$\mathbf{J}_i = R^2(t_1) + z^2(t_1), i = 1, 2.$$

**Example 4.1.** Homing missile guidance with a perfect measurements

$$\tau = 0.01, \kappa_1 = 0.1, \kappa_2 = 0.001, \bar{a}_T^r = 20m/\sec^2, \bar{a}_T^\tau = 20m/\sec^2,$$
$$R(0) = 200m, V_m(0) = 10m/\sec, z(0) = 100, \dot{z}(0) = 40,$$
$$a_T^r(t) = \bar{a}_T^r(\sin(\omega \cdot t))^2, a_T^\tau(t) = \bar{a}_T^\tau \sin(\omega \cdot t), \omega = 5.$$

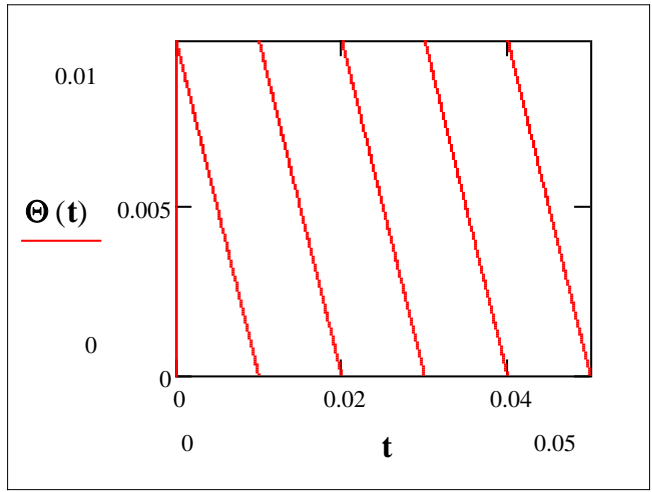

**Pic1.1.** Cutting function $\Theta_\tau(t), \tau = 0.01$.

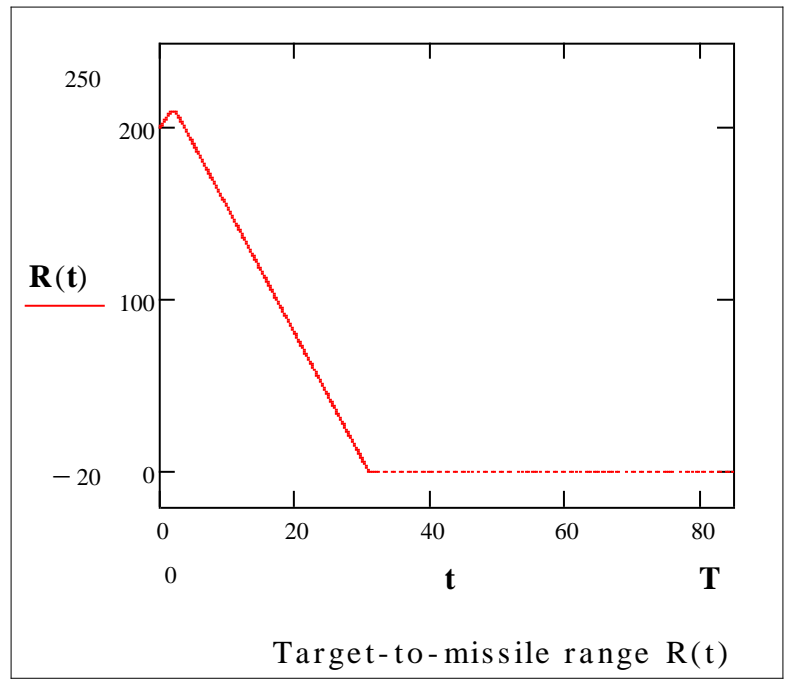

Target-to-missile range R(t)

**Pic1.2.** Target-to-missile range $R(t)$.

$R(T) = 9.179 \times 10^{-8} m.$

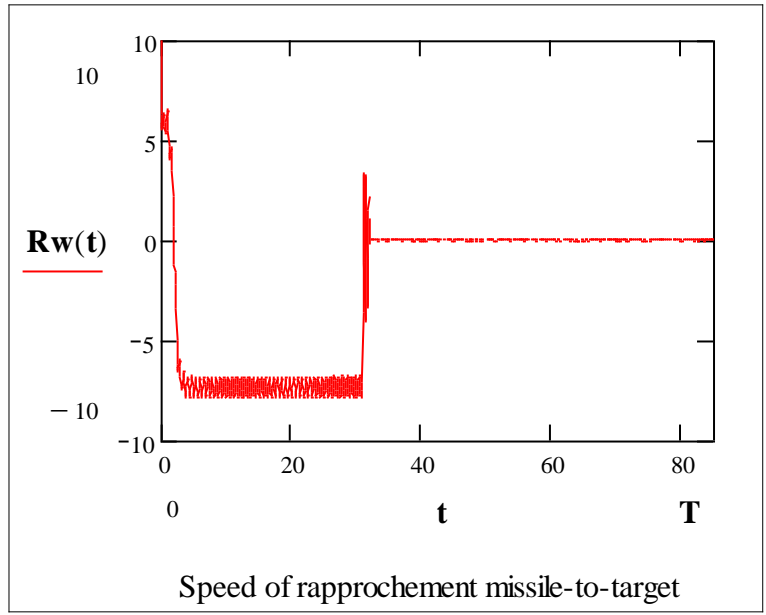

**Pic**.1.3.Speed of rapprochement missile-to-target: $\dot{R}(t)$.
$\dot{R}(0) = 10 m/\sec, \dot{R}(T) = -5.995 \times 10^{-3} m/\sec.$

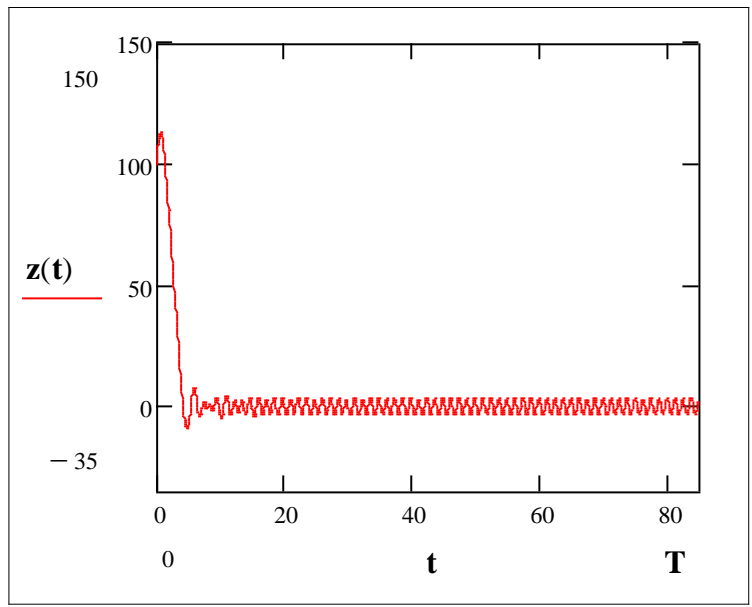

**Pic**1.4.Variable $z(t) = R(t)\dot{\sigma}(t). z(0) = 100, z(T) = 2.922.$

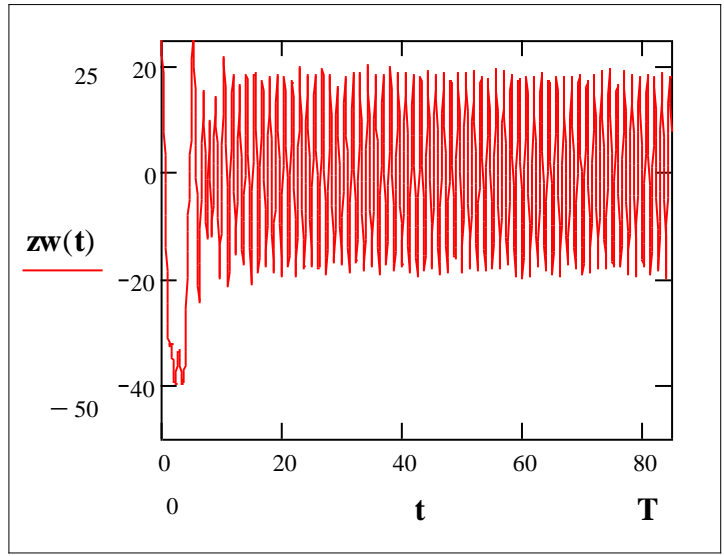

**Pic.3.4**.Variable: $\ddot{z}(t) = \dot{R}\dot{\sigma} + R\ddot{\sigma}.\dot{z}(0) = 40.$

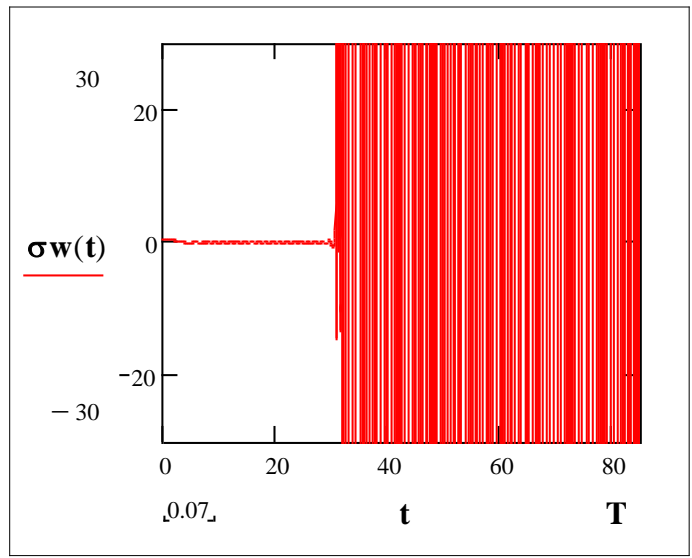

**Pic3.5**.Variable $\dot{\sigma}(t)$.

$\dot{\sigma}(0) = 0.5.$

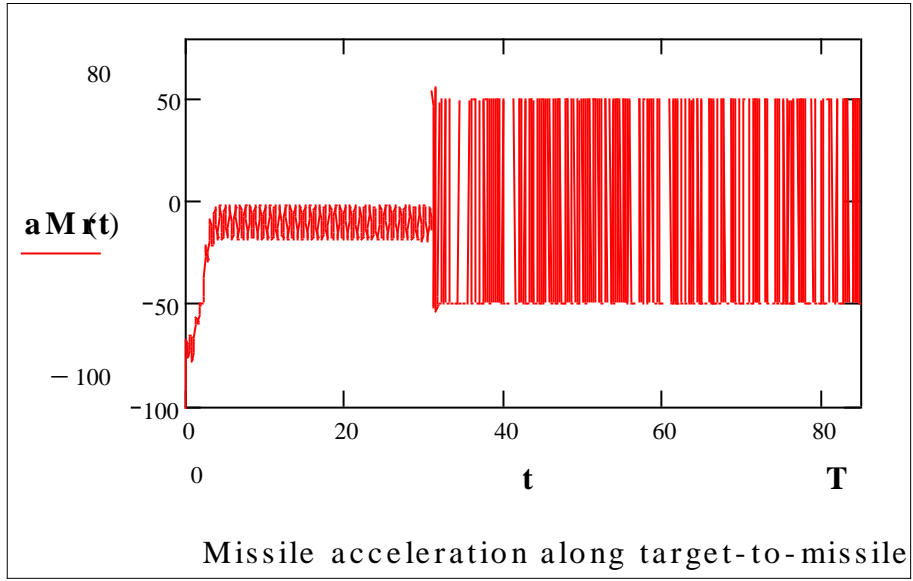

Pic3.6. Missile acceleration along target-to-missile direction $a_M^r(t)$.

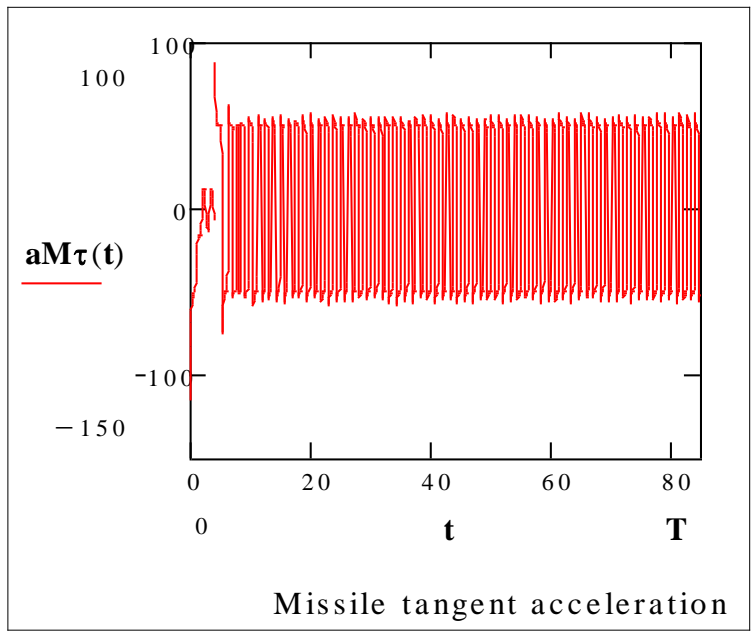

Pic3.7. Missile tangent acceleration $a_M^\tau(t)$.

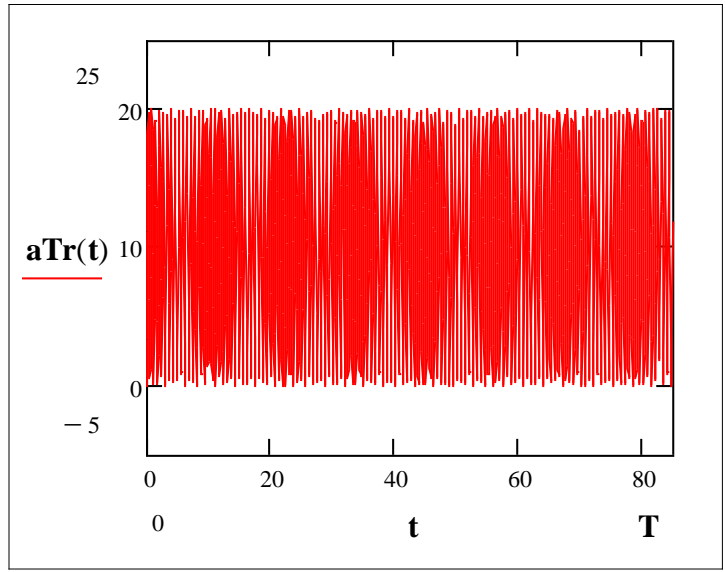

**Pic**.3.8. Target acceleration along target-to-missile direction: $a_T^r(t)$.

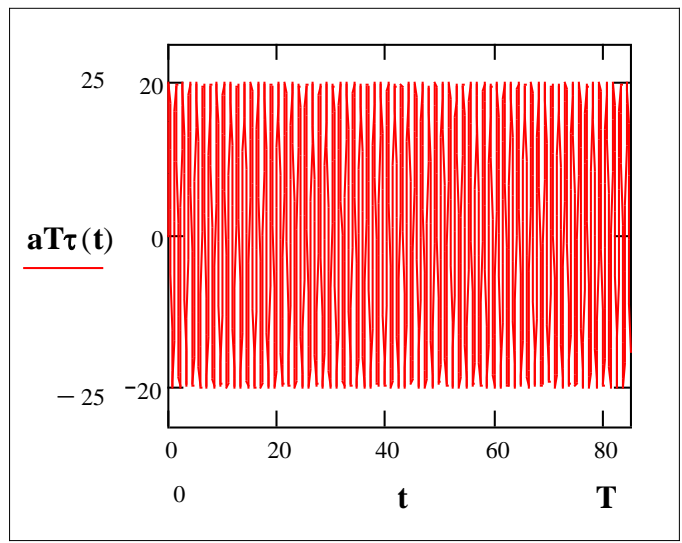

**Pic** 3.9. Target tangent acceleration $a_T^\tau(t)$.

**II. Homing missile guidance with imperfect measurements and imperfect information about the system.**

$\mathbf{f(x)}$ -vector function

$I$ - the variable denotes moment of inertia

$L$ - the variable denotes a lift force

$L_{(\cdot)}$ - the variable denotes a lift force derivative

$M$ - the variable denotes a pitch moment

$M_{(\cdot)}$ - the variable denotes a pitch moment derivative

$m$ - the variable denotes the mass

$q$ - the variable denotes a pitch rate

$\alpha$ - the variable denotes the angle of attack

$\gamma$ - the variable denotes a flight path angle

$\delta$ - the variable denotes the control surface deflection angle

$\boldsymbol{\delta}$ - the variable denotes vector of control surface deflections, canard and tail

$\lambda, \sigma$ - the variables denotes angle between the temporary and initial line of sight

$r, R$ -the variables denotes radial, along the line of sight

$V$ - the variable denotes a speed

$X$ - $M$ - $Z$ -the variable denotes body reference frames

$X$ - $O$ - $Z$ -the variable denotes inertial reference frame

$Z$ -the variable denotes zero-effort miss

$z$ - the variable denotes target-missile relative displacement normal to the initial line of sight

$\delta_c$ - the variable denotes a canard control

$\delta_t$ - the variable denotes a tail control

$\|\mathbf{y}\|$ -denotes euclidian norm of any $\mathbf{y} \in \mathbb{R}^m$

In Fig. 4.2, a schematic view of the planar game geometry is shown.

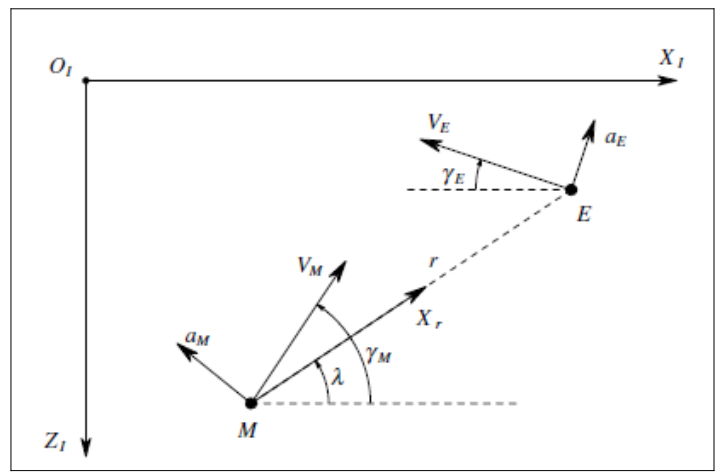

Fig 4.2. Planar game geometry

The engagement kinematics, expressed in a polar coordinate system $(r, \lambda)$ attached to the missile, is

$$\dot{r} = V_r, \dot{\lambda} = V_\lambda/r, \qquad (4.12)$$

where the closing speed $V_r$ is

$$V_r = -[V_M \cos(\gamma_M - \lambda) + V_E \cos(\gamma_M + \lambda)] \qquad (4.13)$$

and the speed perpendicular to the LOS is

$$V_\lambda = -[V_M \sin(\gamma_M - \lambda) - V_E \sin(\gamma_M + \lambda)] \qquad (4.14)$$

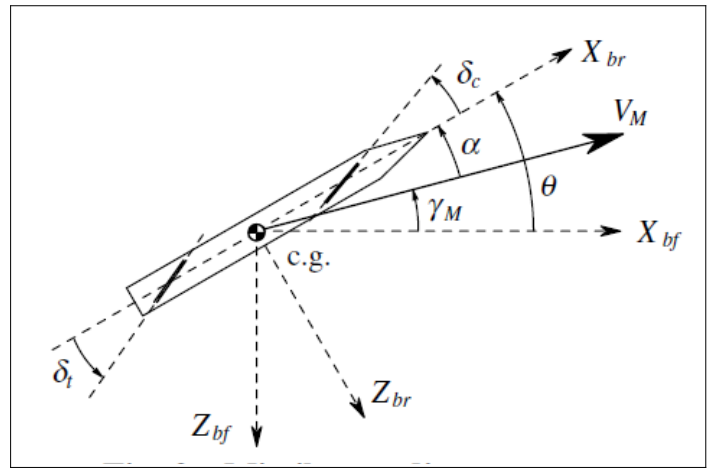

Fig4.3. Missile coordinate systems

The missile planar dynamics are expressed using the coordinate systems presented in Fig.4.3. $X_{bf}$-M-$Z_{bf}$ is parallel to the inertial frame $X_I$--$O_I$--$Z_I$, with its origin located at the missile's center of gravity (c.g.). It is used to express the missile attitude relative to the inertial frame. The missile equations of motion are derived in the rotating body fixed coordinate frame $X_{br}$--M--$Z_{br}$, where the $X_{br}$ axis is aligned with the missile's longitudinal axis. It is assumed that during the end game, the time of interest in our analysis, the missile speed is constant. Thus, the planar missile dynamics in general are given by

$$\dot{\alpha} = q - (L_\epsilon(\alpha,\delta_c,\delta_t))_\epsilon/mV_M + \sqrt{\varepsilon}\dot{W}(t),$$

$$\dot{q} = (M_\epsilon(\alpha,\delta_c,\delta_t))_\epsilon/I + \sqrt{\varepsilon}\dot{W}(t),$$

$$\dot{\theta} = q,$$

$$\dot{\delta}_c = (\delta_c^c - \delta_c)/\tau_c, \dot{\delta}_t = (\delta_t^c - \delta_t)/\tau_t,$$

(4.15)

These surfaces are controlled by actuators, modeled by first-order dynamics with time constants $\tau_c$ and $\tau_t$. The aerodynamic forces and moments are nonlinear, partly unknown functions of the related variables, in particular $\alpha, q, \delta_c$, and $\delta_t$. The missile flight path angle and acceleration perpendicular to the LOS are given by

$$\gamma_M = \theta - \alpha$$

$$a_M^n = a_M \cos(\gamma_M - \lambda)$$

(4.16)

Here, the missile acceleration perpendicular to its velocity vector is given by

$$a_M = (L_\epsilon(\alpha,\delta_c,\delta_t))_\epsilon/m. \tag{4.17}$$

In modeling the missile dynamics, we assume that the lift and the aerodynamic pitch moment in Eqs. (4.12) are generated by the missile body and the control surfaces. This is modeled by

$$L_\epsilon/m = L_\alpha^B g_{1,\epsilon}(\alpha) + L_{\delta_c} \cdot g_{2,\epsilon}(\alpha + \delta_c) + L_{\delta_t} \cdot g_{3,\epsilon}(\alpha + \delta_t),$$

$$M/I = M_\alpha^B \cdot g_{4,\epsilon}(\alpha) + M_q \cdot q + M_{\delta_c} \cdot g_{5,\epsilon}(\alpha + \delta_c) + M_{\delta_t} \cdot g_{6,\epsilon}(\alpha + \delta_t), \tag{4.18}$$

$$L_\alpha^B = L_\alpha - L_{\delta_c} - L_{\delta_t}, M_\alpha^B = M_\alpha - M_{\delta_c} - M_{\delta_t},$$

where $f_i, i = 1\ldots 6$ express the nonlinear aerodynamic characteristics of the missile. We assume now that the true dynamics of the target and missile are unknown to the designer of the missile autopilot and guidance. Thus, only imperfact dynamics can be used, imposing modeling errors. Therefore missile guidance, designed in the sequel, uses an of the nonlinear model with imperfact information about the system, i.e.

$$(\dot{\alpha}_\epsilon)_\epsilon = q - (L_\epsilon(\alpha_\epsilon, \delta_c, \delta_t))_\epsilon / mV_M + (\Delta_{\alpha,\epsilon})_\epsilon + \sqrt{\varepsilon}\dot{W}(t),$$

$$(\dot{q}_\epsilon)_\epsilon = (M_\epsilon(\alpha_\epsilon, q_\epsilon, \delta_c, \delta_t))_\epsilon / I + (\Delta_{q,\epsilon})_\epsilon + \sqrt{\varepsilon}\dot{W}(t), \tag{4.19}$$

$$\dot{\theta} = q, \dot{\delta}_c = (\delta_c^c - \delta_c)/\tau_c,$$

where functions $\Delta_\alpha$ and $\Delta_q$ express the modeling errors, such as unmodeled nonlinearities. The state vector $\mathbf{x}$ of the integrated guidance-autopilot problem is defined by

$$\mathbf{x} = [z, \dot{z}, a_T^n, \alpha, q, \delta_c, \delta_t]^\mathrm{T}, \tag{4.20}$$

where $z$ express target-missile relative displacement normal to the initial line of sight see Fig. 4.3.

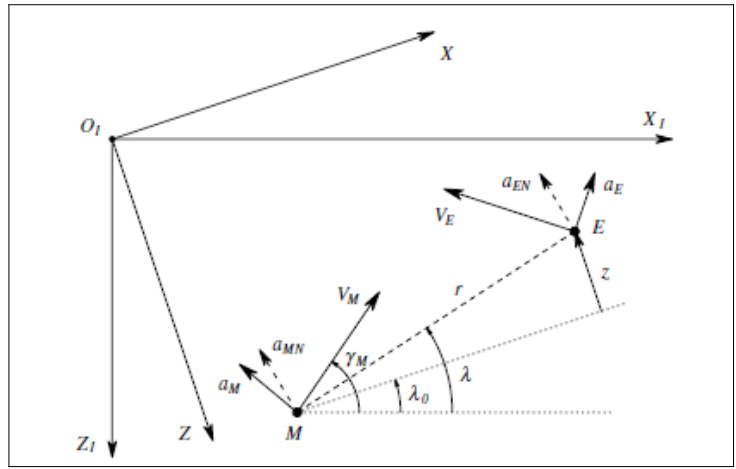

Fig.4.4. Target-missile relative displacement $z$ normal to the initial line of sight

Thus, the missile acceleration normal to the initial LOS is given by

$$a^n_{M_0} = C_{M_0}[\alpha, q, \delta_c, \delta_t]^{\mathrm{T}},$$

$$C_{M_0} = [L_\alpha, 0, L_{\delta_c}, L_{\delta_t}]\cos(\gamma_{M_0} - \lambda_0), \quad (4.21)$$

$$\delta_c = \delta_c(\alpha, z, \dot{z}), \delta_t = \delta_t(\alpha, z, \dot{z})$$

The equations of motion of the integrated dynamics are

$$(\dot{\mathbf{x}}_\epsilon)_\epsilon = A(\mathbf{x}_\epsilon)_\epsilon + (\mathbf{f}_\epsilon(\mathbf{x}_\epsilon))_\epsilon + B\boldsymbol{\delta}^c(\alpha, z, \dot{z}) + Ga^n_T + (\Delta_\epsilon(\mathbf{x}_\epsilon))_\epsilon + \sqrt{\varepsilon}\dot{\mathbf{W}}(t) \quad (4.22)$$

where $\boldsymbol{\delta}^c = [\delta_c, \delta_t]^{\mathrm{T}}$,

$$A = \begin{bmatrix} A_T & A_{12} \\ [0]_{4\times 3} & A_M \end{bmatrix}, A_{12} = \begin{bmatrix} [0]_{1\times 4} \\ -C_{M_0} \\ [0]_{1\times 4} \end{bmatrix}, B = \begin{bmatrix} [0]_{3\times 2} \\ B_M \end{bmatrix}, \quad (4.23)$$

$G = [0, 0, 1/\tau_T, 0, 0, 0, 0]$. Here

$$A_T = \begin{bmatrix} 0 & 1 & 0 \\ 0 & 0 & 1 \\ 0 & 0 & -1/\tau_T \end{bmatrix}, A_M = \begin{bmatrix} -L_\alpha/V_M & 1 & -L_{\delta_c}/V_M & -L_{\delta_t}/V_M \\ M_\alpha & M_q & M_{\delta_c} & M_{\delta_t} \\ 0 & 0 & -1/\tau_c & 0 \\ 0 & 0 & 0 & -1/\tau_t \end{bmatrix},$$

$$(4.24)$$

$$B_M = \begin{bmatrix} 0 & 0 \\ 0 & 0 \\ 1/\tau_c & 0 \\ 0 & 1/\tau_t \end{bmatrix}$$

We assume now that (**i**) CSDE (4.22) is $\widetilde{\mathbb{R}}$-dissipative (see Definition 1) and (**ii**) there exists $r^*$ such that $\forall \mathbf{x}(\|\mathbf{x}\| > r^*)$:

$$\|\Delta_\epsilon(\mathbf{x})\| < O(\|\mathbf{f}_\epsilon(\mathbf{x})\|). \tag{4.25}$$

Taking the origin of the reference frame to be the instantaneous position of the missile, the equation of motion in polar form are:

$$\ddot{R} = R\dot{\sigma}^2 + a_M^r\left[t, \widetilde{\mathbf{R}}(t), \dot{\widetilde{\mathbf{R}}}(t)\right] + a_T^r(t),$$

$$a_M^r(t) \in [-\bar{a}_M^r, \bar{a}_M^r], a_T^r(t) \in [-\bar{a}_T^r, \bar{a}_T^r].$$

$$R\ddot{\sigma} + 2\dot{R}\dot{\sigma} = a_M^n\left[t, \widetilde{\sigma}(t), \dot{\widetilde{\sigma}}(t)\right] + a_T^n(t),$$

$$a_M^n(t) \in [-\bar{a}_M^n, \bar{a}_M^n], a_T^n(t) \in [-\bar{a}_T^n, \bar{a}_T^n] \tag{4.26}$$

1. The variable $R = R(t)$ denotes a true target-to-missle range $R_{TM}(t)$.
2. The variable $\widetilde{\mathbf{R}} = \widetilde{\mathbf{R}}(t)$ denotes the it is *real measured* target-to-missile range $R_{TM}(t)$.
3. The variable $\sigma = \sigma(t)$ denotes a true line-of-sight angle (**LOST**) i.e., the it is *true* angle between the constant reference direction and target-to-missile direction.
4. The variable $\widetilde{\sigma} = \widetilde{\sigma}(t)$ denotes the it is *real measured* line-of-sight angle (**LOSM**) i.e., the it is *true* angle between the constant reference direction and target-to-missile direction.
5. The variable $a_M^n(t) = a_M^n\left[t, \widetilde{\mathbf{R}}(t), \dot{\widetilde{\mathbf{R}}}(t)\right]$ denotes the missiles acceleration along direction which perpendicularly to line-of-sight direction.
6. The variable $a_M^r(t) = a_M^r\left[t, \widetilde{\sigma}(t), \dot{\widetilde{\sigma}}(t)\right]$ denotes the missile acceleration along target-to-missile direction.
7. The variable $a_T^n(t)$ denotes the target acceleration along direction which perpendicularly to line-of-sight direction.
8. The variable $a_T^r(t)$ denotes the target acceleration along target-to-missile direction.

Using replacement $\dot{z} = R\dot{\sigma}$ into Eq.(4.26) one obtain:

$$\ddot{R} = \frac{\dot{z}^2}{R} + a_M^r\left[t, \widetilde{\mathbf{R}}(t), \dot{\widetilde{\mathbf{R}}}(t)\right] + a_T^r(t),$$

$$a_M^r(t) \in [-\bar{a}_M^r, \bar{a}_M^r], a_T^r(t) \in [-\bar{a}_T^r, \bar{a}_T^r].$$

$$\ddot{z} = -\frac{\dot{R}\dot{z}}{R} + a_M^n\left[t, \widetilde{\mathbf{z}}(t), \dot{\widetilde{\mathbf{z}}}(t)\right] + a_T^n(t),$$

(4.27)

$$a_M^n(t) \in [-\bar{a}_M^n, \bar{a}_M^n], a_T^n(t) \in [-\bar{a}_T^n, \bar{a}_T^n].$$

$$\widetilde{\dot{z}}(t) = \widetilde{R}(t)\widetilde{\dot{\sigma}}(t),$$

$$\widetilde{\ddot{z}}(t) = \dot{\widetilde{R}}(t)\widetilde{\dot{\sigma}}(t) + \widetilde{R}(t)\widetilde{\ddot{\sigma}}(t).$$

Suppose that:

$$\widetilde{\mathbf{R}}(t) = R(t) + \beta_1(t),$$

(4.28)

$$\widetilde{\sigma}(t) = \sigma(t) + \beta_2(t).$$

Thus

$$\widetilde{\dot{\mathbf{R}}}(t) = \dot{R}(t) + \dot{\beta}_1(t) = \dot{R}(t) + \bar{\boldsymbol{\beta}}_1(t), \bar{\boldsymbol{\beta}}_1(t) \triangleq \dot{\beta}_1(t).$$

$$\widetilde{\dot{\boldsymbol{\sigma}}}(t) = \dot{\sigma}(t) + \dot{\beta}_2(t), \widetilde{\ddot{\boldsymbol{\sigma}}}(t) = \ddot{\sigma}(t) + \ddot{\beta}_2(t).$$

$$\widetilde{\dot{\mathbf{z}}}(t) = \widetilde{\mathbf{R}}(t)\widetilde{\dot{\boldsymbol{\sigma}}}(t) = (R(t) + \beta_1(t))(\dot{\sigma}(t) + \dot{\beta}_2(t)) =$$

$$= R(t)\dot{\sigma}(t) + [\beta_1(t)(\dot{\sigma}(t) + \dot{\beta}_2(t)) + R(t)\dot{\beta}_2(t)] \approx$$

$$\approx \dot{z}(t) + \left[\beta_1(t)\left(\widetilde{\dot{\boldsymbol{\sigma}}}(t) + \dot{\beta}_2(t)\right) + \widetilde{\mathbf{R}}(t)\dot{\beta}_2(t)\right] =$$

$$= \dot{z}(t) + \bar{\boldsymbol{\beta}}_2(t),$$

$$\bar{\boldsymbol{\beta}}_2(t) \triangleq \beta_1(t)\left(\widetilde{\dot{\boldsymbol{\sigma}}}(t) + \dot{\beta}_2(t)\right) + \widetilde{\mathbf{R}}(t)\dot{\beta}_2(t).$$

(4.29)

$$\widetilde{\ddot{\mathbf{z}}}(t) = \widetilde{\dot{\mathbf{R}}}(t)\widetilde{\dot{\boldsymbol{\sigma}}}(t) + \widetilde{\mathbf{R}}(t)\widetilde{\ddot{\boldsymbol{\sigma}}}(t) =$$

$$(\dot{R}(t) + \dot{\beta}_1(t))(\dot{\sigma}(t) + \dot{\beta}_2(t)) + (R(t) + \beta_1(t))(\ddot{\sigma}(t) + \ddot{\beta}_2(t)) =$$

$$= [\dot{R}(t)\dot{\sigma}(t) + R(t)\ddot{\sigma}(t)] + \dot{\beta}_1(t)(\dot{\sigma}(t) + \dot{\beta}_2(t)) +$$

$$+ \dot{R}(t)\dot{\beta}_2(t) + R(t)\ddot{\beta}_2(t) + \beta_1(t)(\ddot{\sigma}(t) + \ddot{\beta}_2(t)) \approx$$

$$\approx \ddot{z}(t) + \widetilde{\dot{\mathbf{R}}}(t)\dot{\beta}_2(t) + \widetilde{\mathbf{R}}(t)\ddot{\beta}_2(t) + \beta_1(t)\left(\widetilde{\ddot{\boldsymbol{\sigma}}}(t) + \ddot{\beta}_2(t)\right) =$$

$$\ddot{z}(t) + \bar{\boldsymbol{\beta}}_3(t),$$

$$\bar{\boldsymbol{\beta}}_3(t) \triangleq \widetilde{\dot{\mathbf{R}}}(t)\dot{\beta}_2(t) + \widetilde{\mathbf{R}}(t)\ddot{\beta}_2(t) + \beta_1(t)\left(\widetilde{\ddot{\boldsymbol{\sigma}}}(t) + \ddot{\beta}_2(t)\right).$$

Using the replacement $\dot{R} = V_r, \dot{z} = w$ into Eq.(4.27) one obtain:

$$\dot{R} = V_r,$$

$$\dot{V}_r = \frac{\dot{z}^2}{R} + a_M^r\left[t, \widetilde{R}(t), \widetilde{V}_r(t)\right] + a_T^r(t),$$

$$\widetilde{V}_r(t) = \widetilde{\dot{R}}(t) = V_r(t) + \bar{\beta}_1(t),$$

$$a_M^r(t) \in [-\bar{a}_M^r, \bar{a}_M^r], a_T^r(t) \in [-\bar{a}_T^r, \bar{a}_T^r]. \qquad (4.30)$$

$$\dot{z} = w,$$

$$\dot{w} = -\frac{V_r w}{R} + a_M^n\left[t, \widetilde{w}(t), \widetilde{\dot{w}}(t)\right] + a_T^n(t),$$

$$a_M^n(t) \in [-\bar{a}_M^n, \bar{a}_M^n], a_T^n(t) \in [-\bar{a}_T^n, \bar{a}_T^n].$$

Let as considered antagonistic differential game $IDG_{2;T}(f, 0, M, 0, \boldsymbol{\beta}, \mathbf{w})$, $\boldsymbol{\beta} = (\beta_1, \beta_2, \beta_3, \beta_4), \mathbf{w} = (\beta_2, \beta_4)$ with non-linear dynamics and imperfect measurements and imperfect information about the system:

$$\dot{R} = V_r,$$

$$\dot{V}_r = \frac{\dot{z}^2}{R} + a_M^r(t) + a_T^r(t),$$

$$a_M^r(t) = a_M^r\left[t, \widetilde{\mathbf{R}}(t), \widetilde{V}_r(t)\right] - \kappa_1 \widetilde{V}_r^3(t)$$

$$\widetilde{\mathbf{R}}(t) = R(t) + \beta_1(t); \widetilde{V}_r(t) = V_r(t) + \beta_2(t),$$

$$\check{a}_M^r(t) \in [-\bar{a}_M^r, \bar{a}_M^r], a_T^r(t) \in [-\bar{a}_T^r, \bar{a}_T^r].$$

$$\dot{z} = w, \qquad (4.31)$$

$$\dot{w} = -\frac{V_r w}{R} + \check{a}_M^n(t) + a_T^n(t),$$

$$a_M^n(t) = a_M^n\left[t, \widetilde{w}(t), \widetilde{\dot{w}}(t)\right] - \kappa_2(\widetilde{w}(t))^3,$$

$$\widetilde{\mathbf{z}}(t) = \dot{z}(t) + \beta_3(t), \widetilde{\dot{\mathbf{z}}}(t) = \ddot{z}(t) + \beta_4(t),$$

$$a_M^n(t) \in [-\bar{a}_M^n, \bar{a}_M^n], a_T^n(t) \in [-\bar{a}_T^n, \bar{a}_T^n].$$

$$\mathbf{J}_i = R^2(t_1), i = 1, 2.$$

Optimal control problem for the first player are:

$$\bar{\mathbf{J}}_1 = \min_{\check{a}_M^r(t) \in [-\bar{a}_M^r, \bar{a}_M^n], \check{a}_M^\tau(t) \in [-\bar{a}_M^n, \bar{a}_M^n]} \left\{ \max_{a_T^r(t) \in [-\bar{a}_T^r, \bar{a}_T^r], a_T^n(t) \in [-\bar{a}_T^n, \bar{a}_T^n]} R^2(t_1) \right\}. \qquad (4.32)$$

Optimal control problem for the second player are:

$$\bar{\mathbf{J}}_2 = \max_{a_T^r(t) \in [-\bar{a}_T^r, \bar{a}_T^r], a_T^n(t) \in [-\bar{a}_T^n, \bar{a}_T^n]} \left\{ \min_{a_M^r(t) \in [-\bar{a}_M^r, \bar{a}_M^r], a_M^n(t) \in [-\bar{a}_M^n, \bar{a}_M^n]} R^2(t_1) \right\}. \quad (4.33)$$

Let us consider now optimal control problem for the next dissipative nonlinear game with imperfect dinamics:

$$\dot{R} = V_r,$$

$$\dot{V}_r = a_M^r(t) + a_T^r(t),$$

$$a_M^r(t) = a_M^r \left[ t, \widetilde{\mathbf{R}}(t), \widetilde{V}_r(t) \right] - \kappa_1(t) \widetilde{V}_r^3(t),$$

$$\widetilde{\mathbf{R}}(t) = R(t) + \beta_1(t); \widetilde{V}_r(t) = V_r(t) + \beta_2(t),$$

$$a_M^r(t) \in [-\bar{a}_M^r, \bar{a}_M^r], a_T^r(t) \in [-\bar{a}_T^r, \bar{a}_T^r].$$

$$(4.3.34)$$

$$\dot{w} = a_M^n(t) + a_T^n(t),$$

$$a_M^n(t) = a_M^n \left[ t, \widetilde{w}(t), \widetilde{\dot{w}}(t) \right] - \kappa_2 \left( \widetilde{\dot{z}}(t) \right)^3,$$

$$\widetilde{z}(t) = z(t) + \beta_3(t), \widetilde{\dot{z}}(t) = \dot{z}(t) + \beta_4(t)$$

$$a_M^n(t) \in [-\bar{a}_M^n, \bar{a}_M^n], a_T^n(t) \in [-\bar{a}_T^n, \bar{a}_T^n].$$

$$\mathbf{J}_i = R^2(t_1), i = 1, 2.$$

Optimal control problem for the first player are:

$$\bar{J}_1 = \min_{s_M^r(t)\in[-\bar{a}_M^r,\bar{a}_M^n],\, a_M^\tau(t)\in[-\bar{a}_M^n,\bar{a}_M^n]} \left\{ \max_{a_T^r(t)\in[-\bar{a}_T^r,\bar{a}_T^r],\, a_T^n(t)\in[-\bar{a}_T^n,\bar{a}_T^n]} R^2(t_1) \right\}. \tag{4.35}$$

Optimal control problem for the second player are:

$$\bar{J}_2 = \max_{a_T^r(t)\in[-\bar{a}_T^r,\bar{a}_T^r],\, a_T^n(t)\in[-\bar{a}_T^n,\bar{a}_T^n]} \left\{ \min_{a_M^r(t)\in[-\bar{a}_M^r,\bar{a}_M^r],\, a_M^n(t)\in[-\bar{a}_M^n,\bar{a}_M^n]} R^2(t_1) \right\}. \tag{4.36}$$

From Eqs.(4.34)-(4.36) one obtain corresponding linear master game:

$$\dot{r} = v_r + \lambda_2,$$

$$\dot{v}_r = -3\kappa_1 \lambda_2^2 \tilde{v}_r(t) - \kappa_1 \lambda_2^3 + a_M^r(t) + a_T^r(t),$$

$$a_M^r(t) = a_M^r[t, \tilde{r}(t), \tilde{v}_r(t)],$$

$$\tilde{r}(t) = \lambda_1 + r(t) + \beta_1(t); \tilde{v}_r(t) = \lambda_2 + v_r(t) + \beta_2(t),$$

$$a_M^r(t) \in [-\bar{a}_M^r, \bar{a}_M^r], a_T^r(t) \in [-\bar{a}_T^r, \bar{a}_T^r].$$

(4.37)

$$\ddot{z}_1 = -3\kappa_2 \lambda_3^2 \tilde{\dot{z}}_1(t) - \kappa_2 \lambda_3^3 + a_M^n(t) + a_T^n(t),$$

$$a_M^n(t) = a_M^n\left[t, \tilde{z}(t), \tilde{\dot{z}}(t)\right],$$

$$\tilde{z}(t) = \lambda_3 + \dot{z}_1(t) + \beta_3(t), \tilde{\dot{z}}(t) = \dot{z}_1(t) + \beta_4(t)$$

$$a_M^n(t) \in [-\bar{a}_M^n, \bar{a}_M^n], a_T^n(t) \in [-\bar{a}_T^n, \bar{a}_T^n],$$

$$\mathbf{J}_i = r^2(t_1), i = 1, 2.$$

From Eqs.(4.36)-(4.37)
we obtain quasy optimal solution for the antagonistic differential game $IDG_{2;T}(f, 0, M, 0, \boldsymbol{\beta}, \mathbf{w})$ given by Eqs.(4.31-4.34). Quasy optimal control $\{\alpha_M^r(t), \alpha_M^\tau(t)\}$ of the first player and quasy optimal control $\{\alpha_T^r(t), \alpha_T^\tau(t)\}$ of the second player are:

$$\check{\alpha}_M^r(t) \simeq -\rho_M^r \mathbf{sign}[[R(t) + \beta_1(t)] + \Theta_\tau(t)[V_r(t) + \beta_2(t)]] - \kappa_1[V_r(t) + \beta_2(t)]^3,$$

$$\check{\alpha}_M^n(t) = -\rho_M^n \mathbf{sign}[[z(t) + \beta_3(t)] + \Theta_\tau(t)[\dot{z}(t) + \beta_4(t)]] - \kappa_2[\dot{z}(t) + \beta_4(t)]^3.$$

$$\alpha_T^r(t) \simeq \rho_T^r \mathbf{sign}\big[[R(t) + \hat{\beta}_1(t)] + \Theta_\tau(t)[V_r(t) + \hat{\beta}_2(t)]\big] - \kappa_1\big[V_r(t) + \hat{\beta}_2(t)\big]^3, \qquad (4.38)$$

$$\alpha_T^n(t) = \rho_T^n \mathbf{sign}\big[[z(t) + \hat{\beta}_3(t)] + \Theta_\tau(t)[\dot{z}(t) + \hat{\beta}_4(t)]\big] - \kappa_2\big[\dot{z}(t) + \hat{\beta}_4(t)\big]^3.$$

Thus for numerical simulation we obtain ODE:

$$\dot{R} = V_r,$$

$$\dot{V}_r = \frac{\dot{z}^2}{R} - \rho_M^r \mathbf{sign}[[R(t) + \beta_1(t)] + \Theta_\tau(t)[V_r(t) + \beta_2(t)]] -$$
$$-\kappa_1[V_r(t) + \beta_2(t)]^3 + a_T^r(t), \qquad (4.39)$$

$$\ddot{z} = -\frac{V_r \dot{z}}{R} - \rho_M^n \mathbf{sign}[[z(t) + \beta_3(t)] + \Theta_\tau(t)[\dot{z}(t) + \beta_4(t)]] -$$
$$-\kappa_2[\dot{z}(t) + \beta_4(t)]^3 + a_T^n(t).$$

**Example 4.** $\tau = 0.001, \kappa_1 = 10^{-3}, \kappa_2 = 0.001, \bar{a}_T^r = 20 m/\sec^2,$
$\bar{a}_T^\tau = 20 m/\sec^2, \quad R(0) = 200m, V_r(0) = 10 m/sec,$
$z(0) = 60, \dot{z}(0) = 40, a_T^r(t) = \bar{a}_T^r (\sin(\omega \cdot t))^p,$
$a_T^\tau(t) = \bar{a}_T^\tau (\sin(\omega \cdot t))^q, \omega = 50, w(t) = \beta(t) = \bar{\beta}(\sin(\omega \cdot t))^p,$

$\bar{\beta} = 20, \, p = 2, q = 1.$

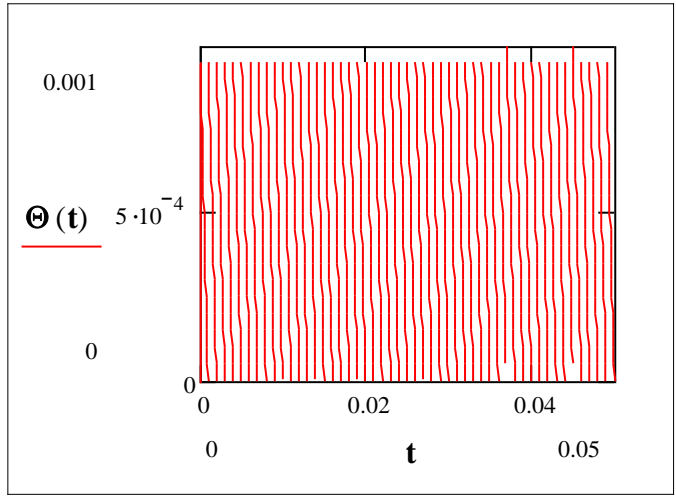

**Pic4.1.** Cutting function $\Theta_\tau(t)$.

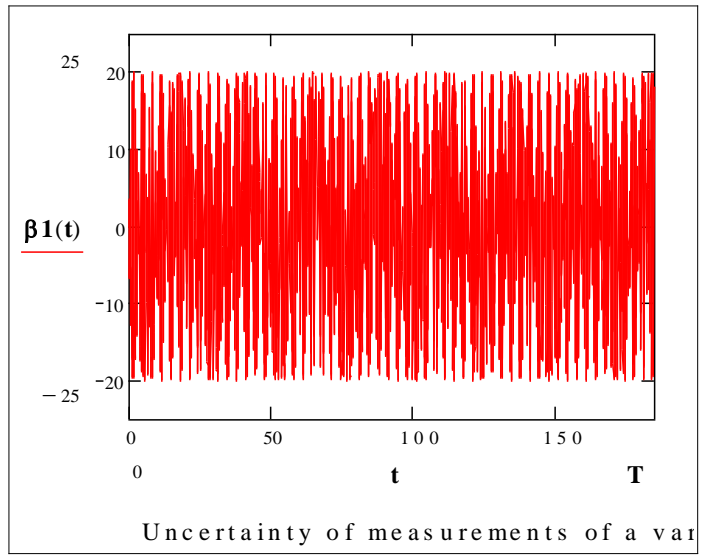

**Pic4.2.** Uncertainty of measurements of a variable $\dot{e}(t)$: $\beta_1(t)$.

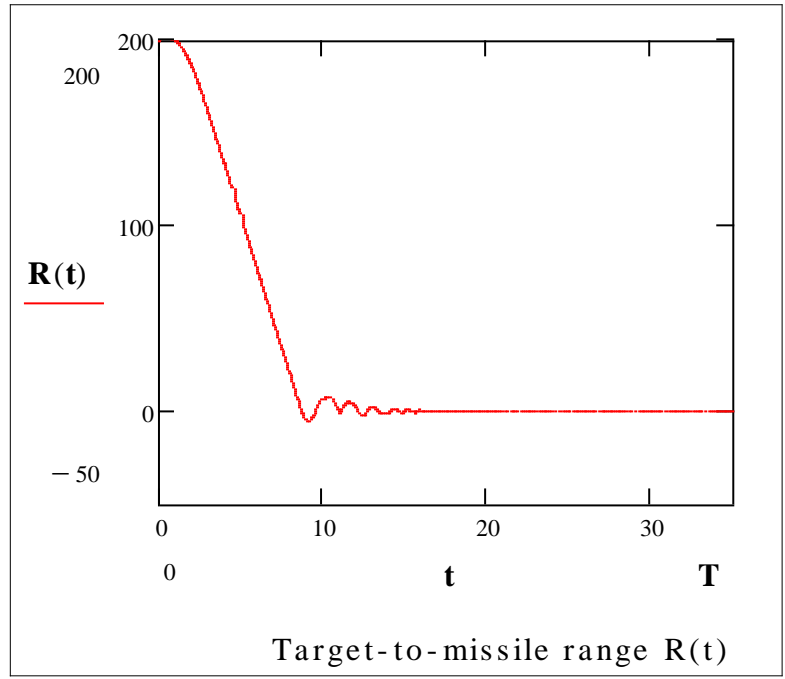

**Pic4.3**. Target-to-missile range $R(t)$.

$$R(T) = 7.2 \times 10^{-3} m.$$

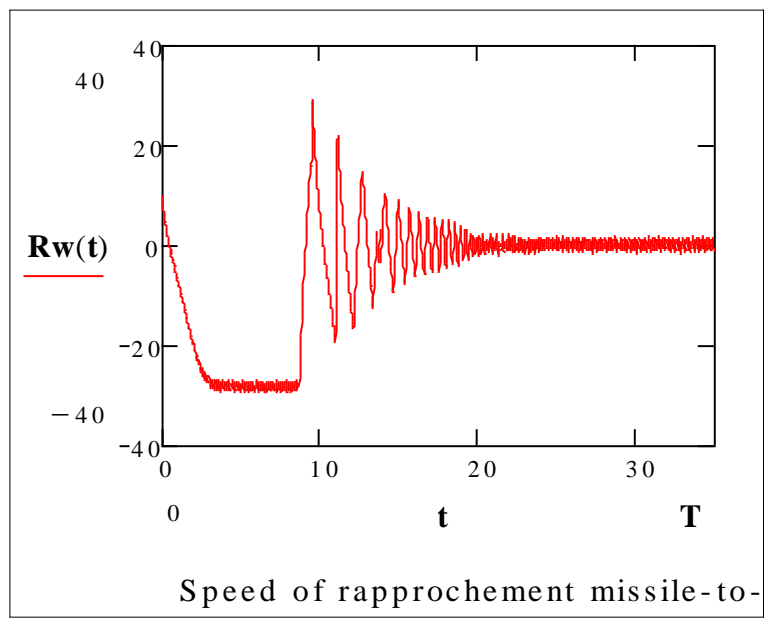

Speed of rapprochement missile-to-

**Pic4.4**. Speed of rapprochement missile-to-target $\dot{R}(t)$.

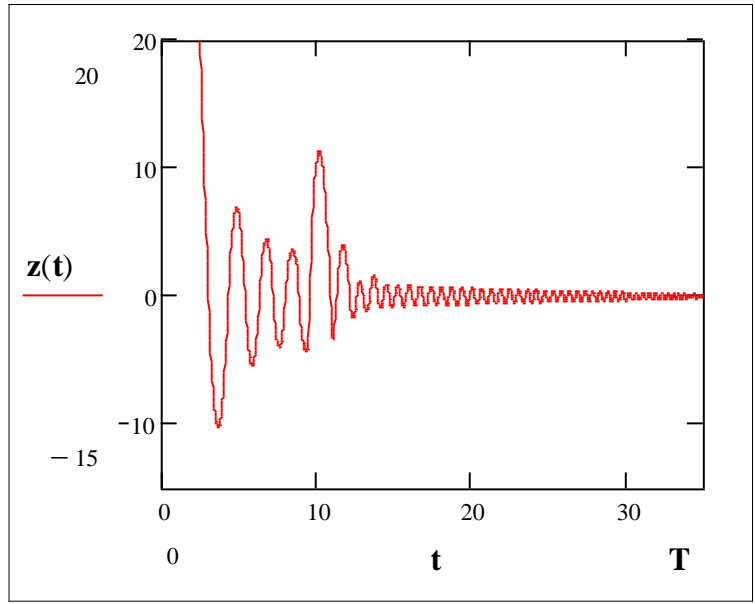

**Pic4.5**. Variable $z(t)$.

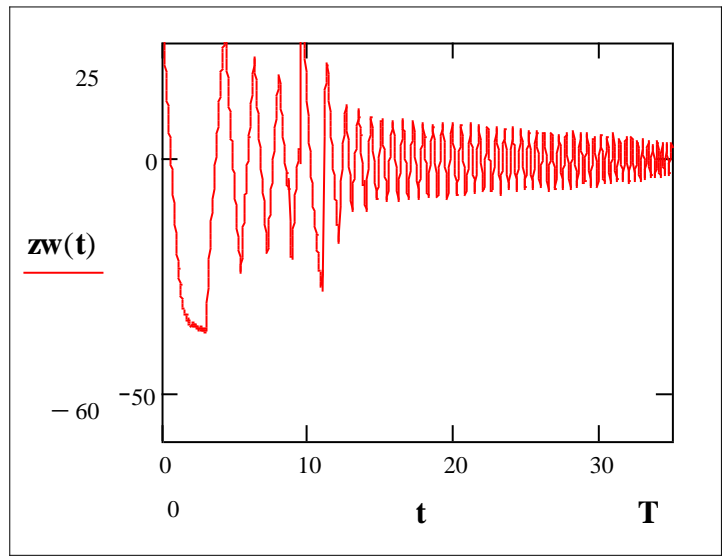

**Pic4.6**. Variable $\dot{z}(t)$. $\dot{z}(T) = 2.172$.

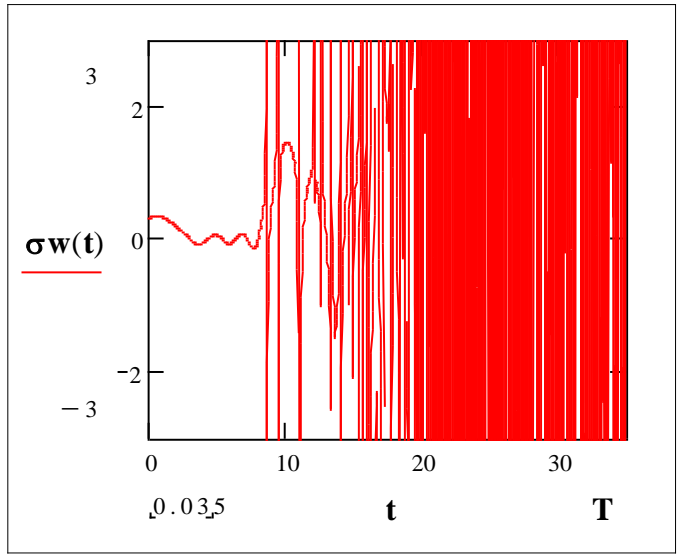

**Pic 4.7.** Variable $\dot{\sigma}(t)$. $\dot{\sigma}(0) = 0.3$.

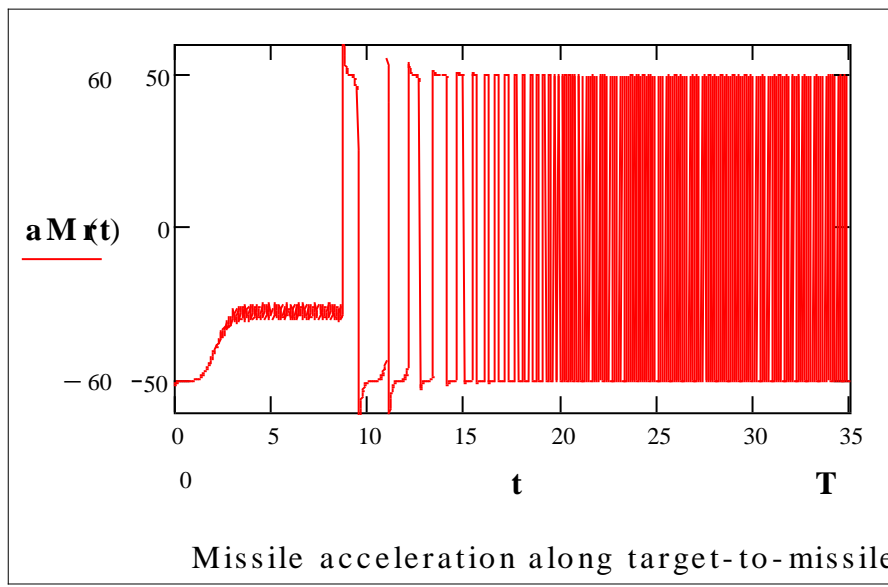

Missile acceleration along target-to-missile

**Pic 4.8.** Missile acceleration along target-to-missile direction $a_M(t)$.

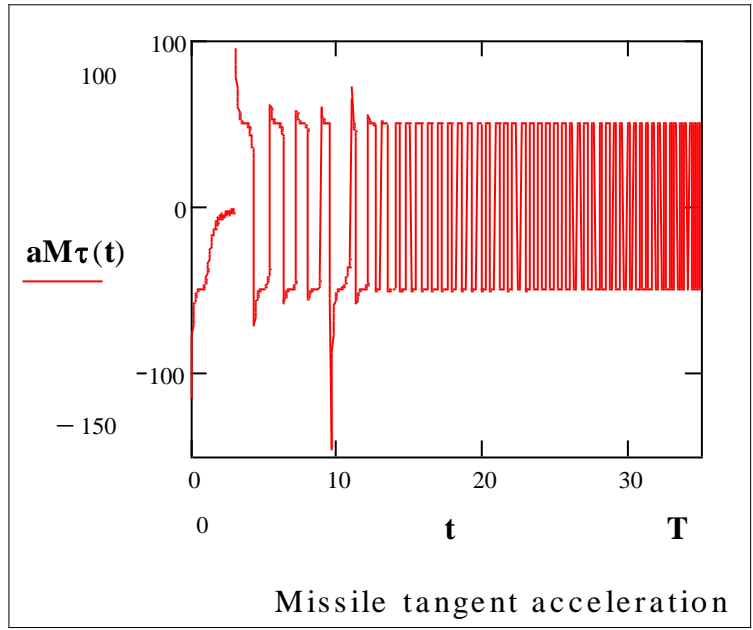

**Pic4.9.** Missile acceleration along direction which perpendicular to line-of-sight direction $a_M^\eta(t)$.

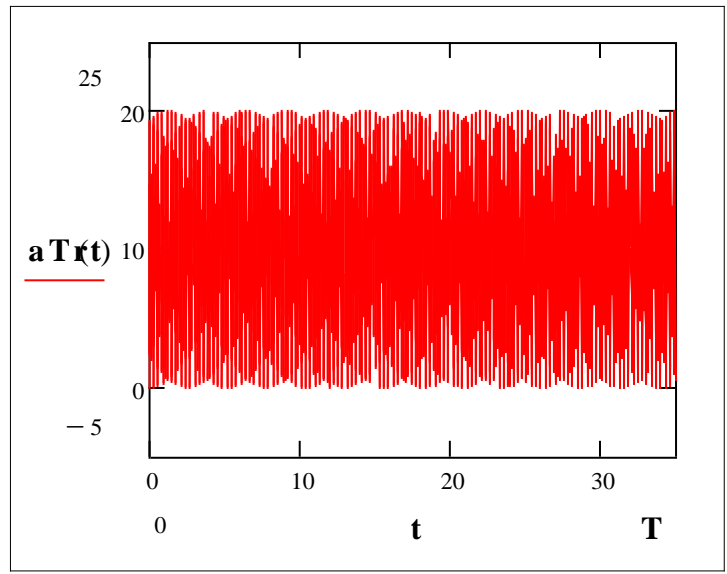

**Pic4.10.** Target acceleration along target-to-missile direction $a_T^r(t)$.

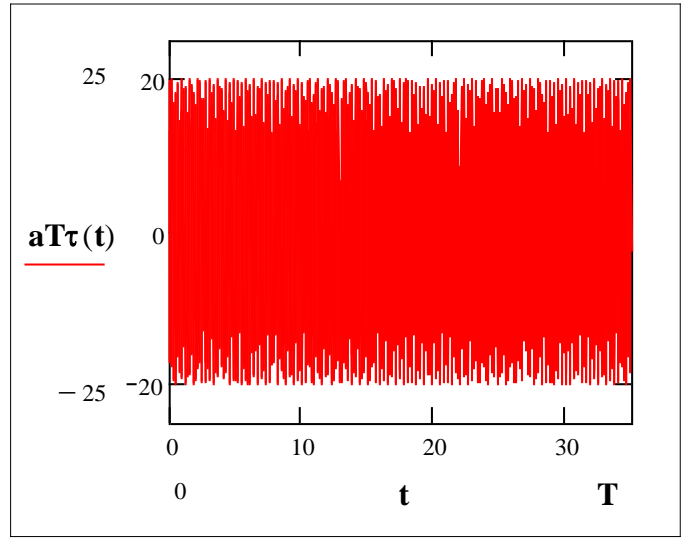

**Pic4.11.** Target acceleration along direction which perpendicular to lyne of-sight direction $a_T^\tau(t)$.

**Example 5.** $\tau = 0.001, \kappa_1 = 10^{-4}, \kappa_2 = 10^{-3}, \bar{a}_T^r = 20 m/\sec^2$,
$\bar{a}_T^\tau = 20 m/\sec^2, \quad R(0) = 200 m, V_r(0) = 10 m/\sec$,
$z(0) = 60, \dot{z}(0) = 40, a_T^r(t) = \bar{a}_T^r (\sin(\omega \cdot t))^p$,
$a_T^\tau(t) = \bar{a}_T^\tau (\sin(\omega \cdot t))^q, \omega = 50, w(t) = \beta_1(t) = \bar{\beta}_1 (\sin(\omega \cdot t))^p$,

$\bar{\beta}_1 = 200 m/\sec, p = 2, q = 1.$

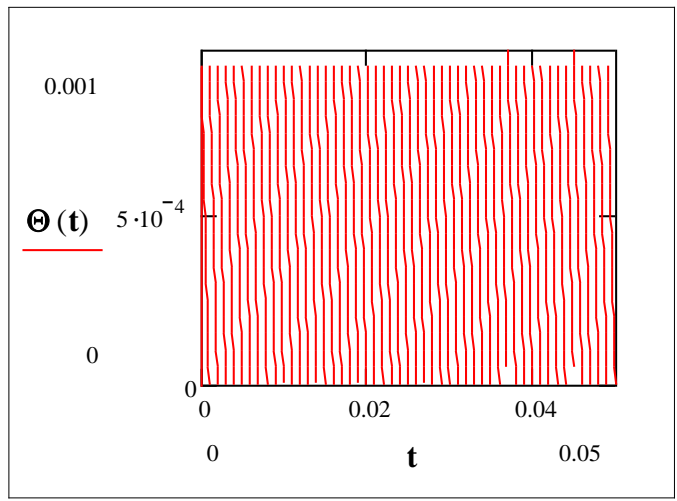

**Pic5.1.** Cutting function $\Theta_\tau(t)$.

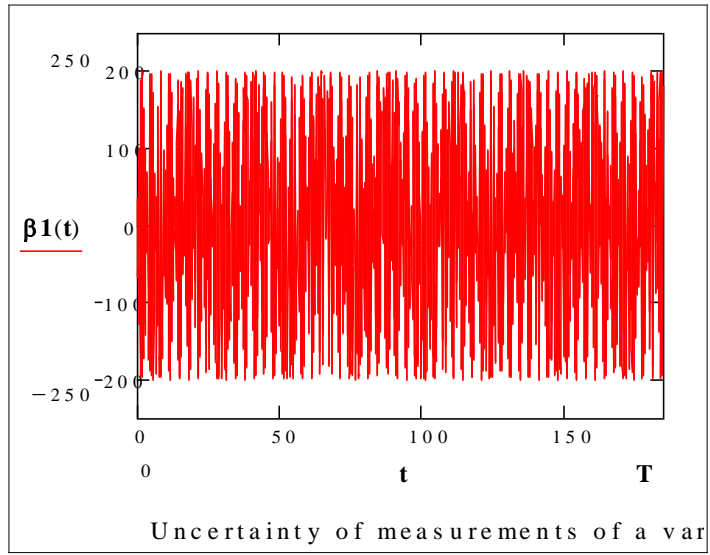

**Pic5.2**.

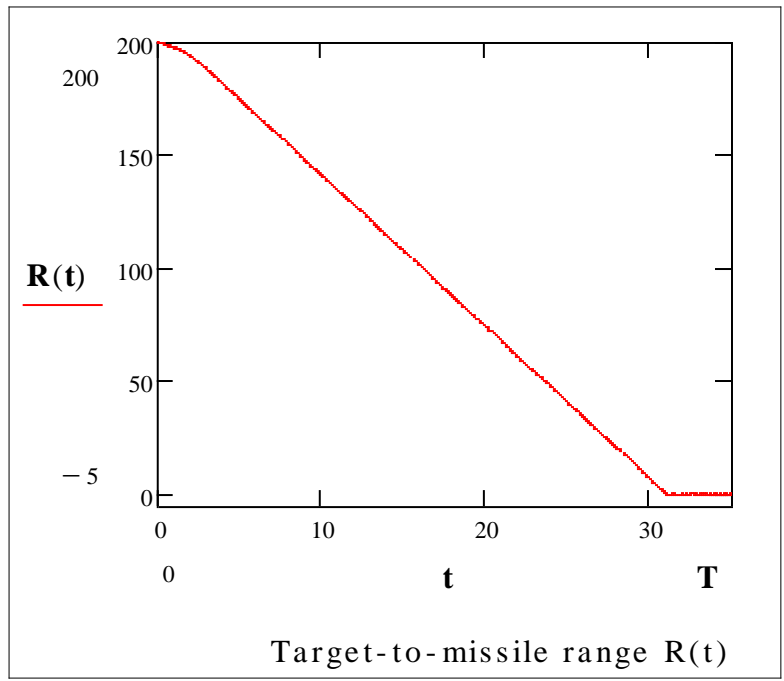

**Pic5.3**. Target-to-missile range $R(t)$. $R(T) = 0.055m$.

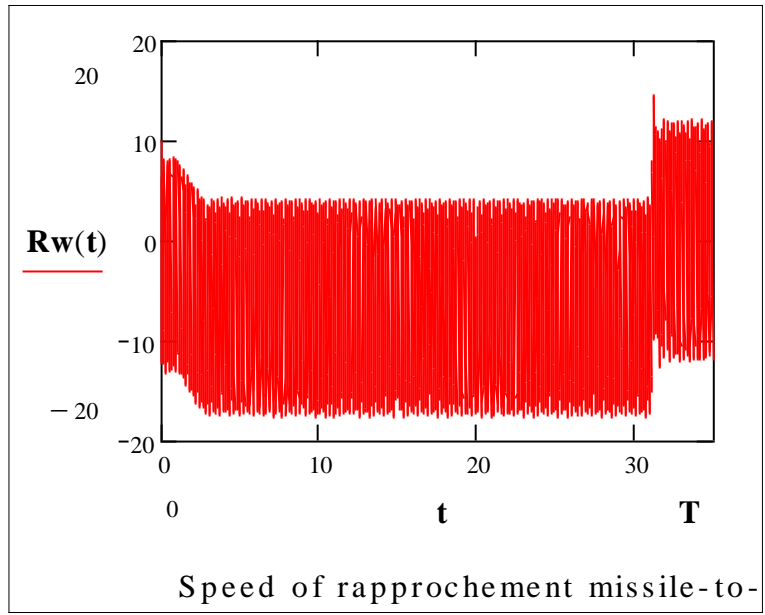

**Pic 5.4.** Speed of rapprochement missile-to-target $\dot{R}(t)$.

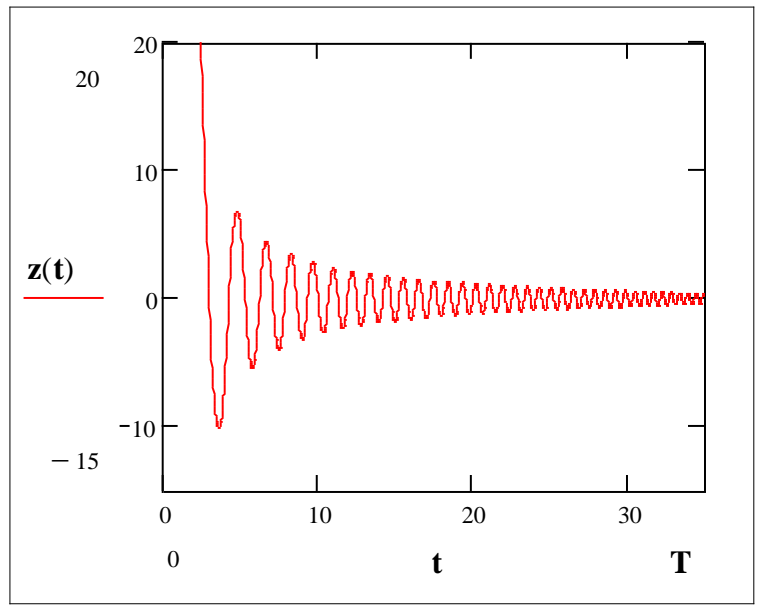

**Pic 5.5.** Variable $z(t)$. $\dot{z}(T) = 0.42$.

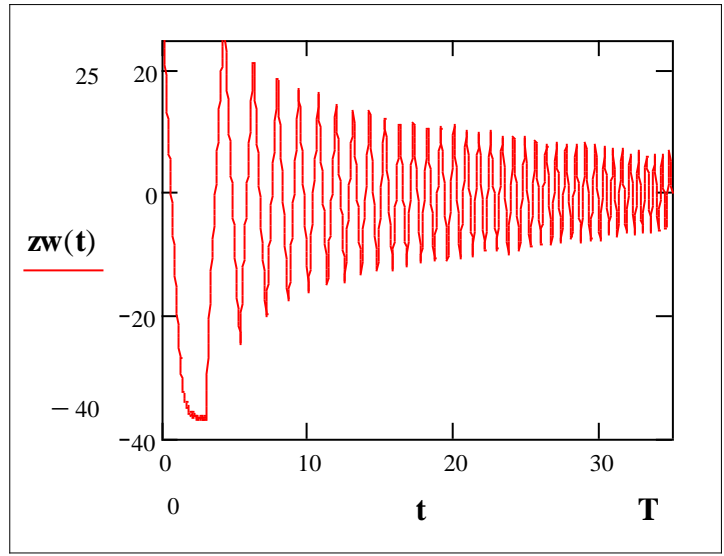

**Pic 5.6.** Variable $z(t)$. $\ddot{z}(T) = -0.149$.

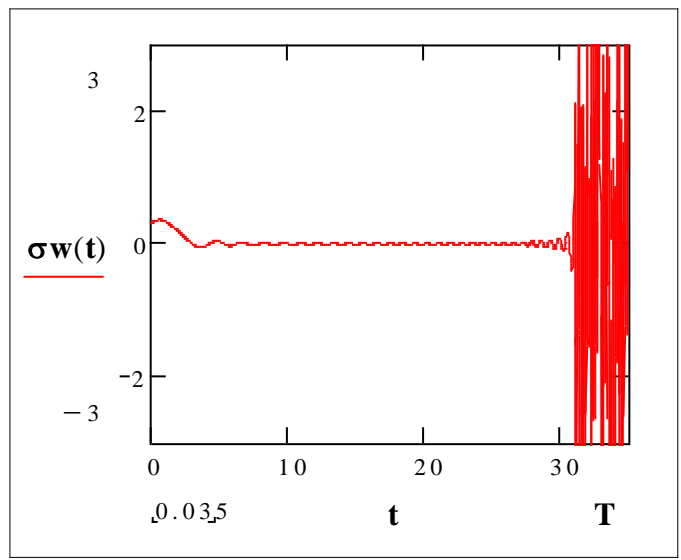

**Pic 5.7.** Variable $\dot{e}(t)$.

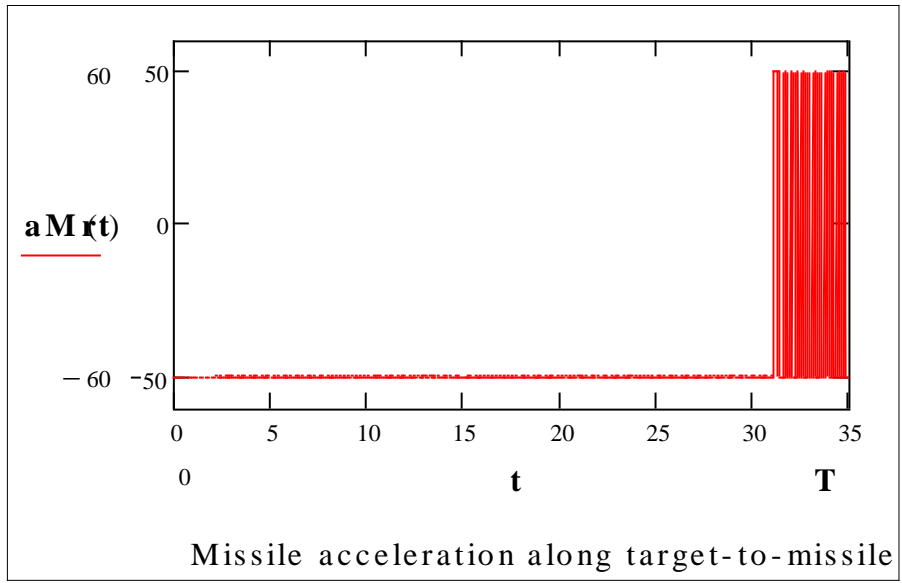

**Pic2.8.** Missile acceleration along target-to-missile direction $a^r_M(t)$.

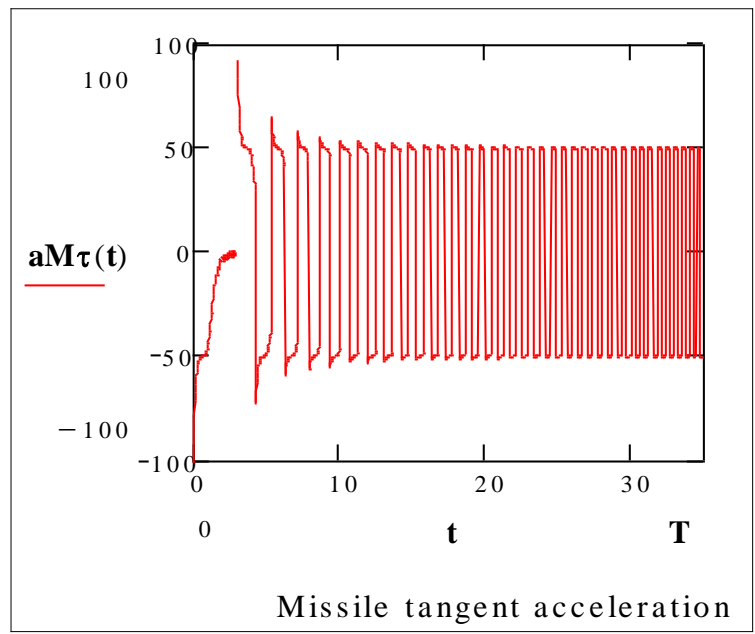

**Pic5.9.** Missile acceleration along direction which perpendicular to line-of-sight direction $a^\tau_M(t)$.

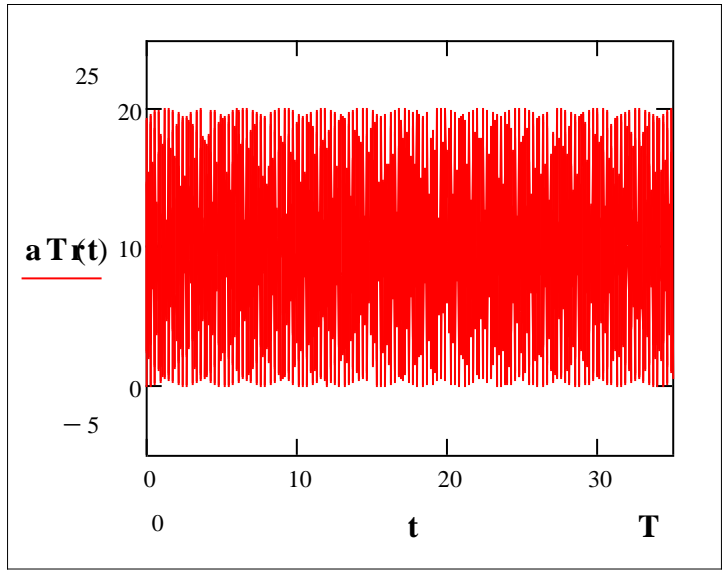

**Pic 5.10.** Target acceleration along target-to-missile direction $a_T^r(t)$.

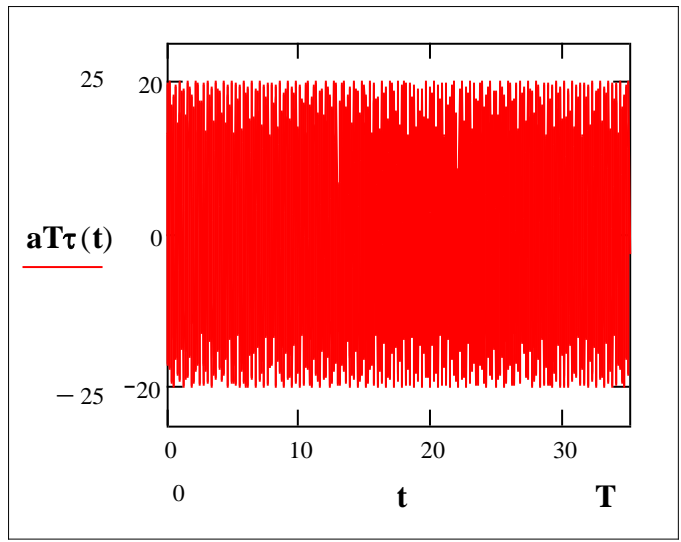

**Pic 5.11.** Target acceleration along direction which perpendicular to line-of-sight direction $a_T^\tau(t)$.

**References.**